\def\bibart#1#2#3#4#5#6#7
{\bibitem{#1}
#2, #4, #5 #6 (#3), #7}

\def\bibcoll#1#2#3#4#5#6#7#8
{\bibitem{#1}
#2, #4, in: #5, #6, #7, #3, pp. #8}

\def\bibbook#1#2#3#4#5#6
{\bibitem{#1} #2, #4, #5, #6, #3}

\def\bibdiss#1#2#3#4#5#6
{\bibitem{#1}
#2 (#3), #4, #5, #6}

\def\updownharpoons{\upharpoonleft \! \downharpoonright }

\def\hochdamit{\! \Rsh }

\def\mat#1{\mathbb{#1}}

\def\qed{\hfill\ \rule{2mm}{2mm} }

\def\qex{\hfill\ \vbox{\hrule\hbox{\vrule\kern4pt\vbox{\kern4pt{}
\kern4pt}\kern4pt\vrule}\hrule}}

\newcounter{casectr}

\newcounter{reductctr}

\newcounter{claimctr}

\newcounter{classctr}

\documentclass[12pt]{article}

\usepackage{amssymb}

\newtheorem{guess}{Guess}[section]

\newtheorem{define}[guess]{Definition}
\newtheorem{prop}[guess]{Proposition}
\newtheorem{theorem}[guess]{Theorem}

\newtheorem{cor}[guess]{Corollary}

\newtheorem{lem}[guess]{Lemma}

\newtheorem{nota}[guess]{Notation}

\newtheorem{remark}[guess]{Remark}

\newtheorem{exam}[guess]{Example}

\newtheorem{oques}[guess]{Open Question}

\usepackage[pdftitle={Small Width Automorphism Conjecture}, pdfauthor={Bernd Schroeder},
colorlinks=true, pdfstartview=FitV, linkcolor=blue, citecolor=blue, urlcolor=blue]{hyperref}

\usepackage{color}

\begin{document}

\bibliographystyle{plain}

\title{The Automorphism Conjecture for
Ordered Sets of Width $\leq 11$}

\author{
\small Bernd S. W. Schr\"oder \\
\small School of Mathematics and Natural Sciences\\
\small The University of Southern Mississippi\\
\small
118 College Avenue, \#5043\\
\small Hattiesburg, MS 39406\\
}

\date{\small \today}

\maketitle

\begin{abstract}

We introduce a recursive method to deconstruct
the automorphism group of an ordered set.
By connecting this method with
deep results for permutation groups,
we prove the Automorphism Conjecture
for ordered sets of width less than or equal to $11$.
Subsequent investigations show that the method presented
here
could lead to a resolution of the Automorphism Conjecture.

\end{abstract}


\section{Introduction}

An {\bf ordered set}
consists of an underlying set $P$ equipped with a
reflexive, antisymmetric and transitive relation $\leq $, the order
relation.
An {\bf order-preserving self-map},
or, an {\bf endomorphism},
of an ordered set $P$ is a
self map
$f:P\to P$ such that $p\leq q$ implies
$f(p)\leq f(q)$.
Consistent with standard terminology,
endomorphisms with an inverse that is
an endomorphism, too, are called {\bf automorphisms}.
The set of endomorphisms is denoted ${\rm End} (P)$
and the set of automorphisms is denoted ${\rm Aut} (P)$.
Rival and Rutkowski's
Automorphism Problem (see \cite{RiRut}, Problem 3)
asks the following.

\begin{oques}

{\bf (Automorphism Problem.)}
Is it true that
$$\lim _{n\to \infty } \max _{|P|=n}
{ |{\rm Aut} (P)|\over |{\rm End} (P)| } =0?$$

\end{oques}

The {\bf Automorphism Conjecture} states that
the
Automorphism Problem
has an affirmative answer.
In light of the facts that, for
almost every ordered set, the identity is the only
automorphism
(see \cite{Proe}, Corollary 2.3a), and that
every ordered set has at least $2^{{2\over3} n} $ endomorphisms
(see \cite{DRSW}, Theorem 1), this conjecture
is quite natural.
Indeed, if, for ordered sets with ``many" automorphisms,
we could show that
there are
``enough" endomorphisms to guarantee the ratio's convergence to zero
(for examples of this technique, see \cite{LRZ,LiuWan}, or
Proposition \ref{manyinmaxlock} here), the conjecture
would be confirmed.
However, the Automorphism Conjecture has been remarkably resilient
against attempts to prove it in general.

Recall that an {\bf antichain} is an ordered set in which no two elements are
comparable and that the {\bf width} $w(P)$
of an ordered set $P$
is the size of the largest antichain contained in $P$.
The Automorphism Conjecture for ordered sets of small width
has recently gathered attention in \cite{BonaMartin}.
It is easy to slightly improve
Theorem 1 in \cite{DRSW} to,
for ordered sets of bounded width,
provide at least $2^{(1-\varepsilon ) n} $ endomorphisms,
see Lemma \ref{2tonforbddwidth}.
With such lower bounds available,
it is natural to also consider upper bounds on the number of
automorphisms.
We will see here that the search for
upper bounds on the number of automorphisms
is linked with
numerous insights on
the connection between the combinatorial
structure of an ordered set
and the structure of its automorphism group.

We start our investigation with
ordered sets that have a lot of local symmetry
in Section \ref{maxlocksec}.
Proposition \ref{maxlockdef}
essentially shows that, if too much local symmetry
is allowed, then, for any automorphism,
the remainder of the ordered set is locked into
following the automorphism's action on a small subset.
Section \ref{orbclustsec} provides an overall framework
to investigate
this ``transmission of local actions of automorphisms."
Proposition \ref{getallfromiou}
shows that the framework decomposes into subsets,
called interdependent orbit unions,
on which the actions of automorphisms
are independent
of each other.
Section \ref{orbgraphsec} focuses on the automorphism groups of
such interdependent orbit unions.
Theorem \ref{pruneorbit6} is the key to splitting the
automorphism group into two parts, which can then be
analyzed separately.
Via the deconstruction of automorphisms
in Section \ref{AuttoLambda},
we can leverage the substantial body of work on
permutation groups.
In Section \ref{primnestsec},
Theorem \ref{permgrtoordaut} then uses bounds for
permutation groups induced on individual orbits to
obtain a bound for $|{\rm Aut} (P)|$ for many types of ordered sets.
Section \ref{ACsmwidthproof}
shows the utility of this approach by, in Theorem \ref{ACw<=12},
confirming the
Automorphism Conjecture for ordered sets of width $\leq 11$.
Section \ref{moreendosec} improves the
lower bound for the number of endomorphisms from \cite{DRSW}
in Theorem \ref{moreendos}
to $2^{0.7924n} $ for large $n$.
Finally,
Section \ref{primorbsec}
gives further upper bounds on the number of automorphisms
for some classes of ordered sets.
The overall total of these advances suggests that
a resolution of the Automorphism Conjecture may well be
within our reach.

\section{Max-locked Ordered Sets}
\label{maxlocksec}

An ordered set with ``many" automorphisms must have a high degree of symmetry.
In terms of counting techniques, this means that, even when
the values of an automorphism are known for a ``large" number of points,
the automorphism would still not be uniquely determined by these values.
It is thus natural to try to identify
smaller subsets $S\subset P$
such that
every automorphism
of the ordered set $P$ is uniquely determined
by its values on $S$.

For sets $B,T\subseteq P$, we will write
$B<T$ iff every $b\in B$ is strictly below every $t\in T$.
For singleton sets, we will omit the set braces.
Recall that a nonempty subset $A\subseteq P$ is called
{\bf order-autonomous}
iff, for all $z\in P\setminus A$, we have that
existence of an $a\in A$ with $z<a$ implies $z<A$, and,
existence of an $a\in A$ with $z>a$ implies $z>A$.
An
order-autonomous subset $A\subseteq P$ will be called {\bf nontrivial}
iff $|A|\not\in \{ 1,|P|\} $.

Lemma \ref{allbutoneac} below shows the
simple idea that we will frequently use:
When there are no nontrivial order-autonomous antichains, then
an automorphism's values
on an antichain are determined by
the
automorphism's
values away from the antichain.
Recall that
$\uparrow x=\{ p\in P:p\geq x\} $
and
$\downarrow x=\{ p\in P:p\leq x\} $.

\begin{lem}
\label{allbutoneac}

Let $P$ be an ordered set with at least 3 points,
let
$A\subseteq P$ be an antichain
that does not contain any nontrivial order-autonomous antichains,
and let
$\Phi , \Psi \in {\rm Aut}
(P)$
be so that both
$\Phi $ and $\Psi $ map $A$ to itself and
$\Phi |_{P\setminus A} =
\Psi |_{P\setminus A} $.
Then $\Phi =\Psi $.

\end{lem}

{\bf Proof.}
Let
$\Phi , \Psi \in {\rm Aut}
(P)$
be so that both
$\Phi $ and $\Psi $ map $A$ to itself and
$\Phi |_{P\setminus A} =
\Psi |_{P\setminus A} $, and let
$x\in A$.
Then
$\uparrow \Phi (x)\setminus \{ \Phi (x)\}
=
\Phi [\uparrow x\setminus \{ x\} ]
=
\Phi |_{P\setminus A} [\uparrow x\setminus \{ x\} ]
=
\Psi |_{P\setminus A} [\uparrow x\setminus \{ x\} ]
=
\uparrow \Psi (x)\setminus \{ \Psi (x)\}
$, and similarly,
$\downarrow \Phi (x)\setminus \{ \Phi (x)\}
=
\downarrow \Psi (x)\setminus \{ \Psi (x)\}
$.
Hence,
$\{ \Phi (x), \Psi (x)\} \subseteq A$
is an order-autonomous antichain.
By hypothesis,
$\{ \Phi (x), \Psi (x)\} \subseteq A$
is not a nontrivial order-autonomous antichain.
Because
there is at least one more point, we
have $\{ \Phi (x), \Psi (x)\} \not= P$ and hence
$\Phi (x)= \Psi (x)$.
\qed

\vspace{.1in}

Recall that an element $x$ of a finite ordered set $P$ is
said to be {\bf minimal} or of {\bf rank 0}, and we set
${\rm rank} (x):=0$,
iff there is no $z\in P$ such that $z<x$.
Recursively, the element $x$ is said to be of
{\bf rank $k$}, and we set
${\rm rank} (x):=k$,
iff $x$ is minimal in $P\setminus \{ z\in P: {\rm rank} (z)\leq k-1\} $.
It is easy to see that the rank
of a point is preserved by automorphisms.

\begin{define}
\label{inducedonrankdef}

Let $P$ be
an
ordered set.
For every nonnegative integer $j$, we define
$R_j $ to be the set of elements of rank $j$.
The largest number $h$ such that $R_h \not= \emptyset $
is called the {\bf height} of $P$.

\end{define}

Clearly,
the symmetric group $S_n $ and the alternating group
$A_n $ are the largest permutation groups
on $n$ elements.
In fact, the very deep
Corollary 1.4 in \cite{Maroti}
(see Theorem \ref{24special} here) shows that these two are
{\em by far}
the largest permutation groups
on $n$ elements.
Hence we start our analysis with ordered sets on which
a set $R_k $ carries a
maximum-sized alternating group.
We will return to a more detailed investigation of
permutation groups
on antichains
in Section \ref{primnestsec}.

\begin{define}

For any set $S$, we define $A_{|S|} (S)$ to be the
alternating group on $S$.

\end{define}

\begin{lem}
\label{noverklem}

(Folklore.)
Let $n\geq 6$ be a composite number,
let $k$ be a nontrivial divisor of $n$ and
let $\ell \leq k$ be the smallest nontrivial
divisor of $n$.
Then
$
k!\left( \left( {n\over k} \right) ! \right) ^k
\leq
\ell !\left( \left( {n\over \ell } \right) ! \right) ^\ell
< (n-1)!$.

\end{lem}

{\bf Proof.}
For the first inequality,
every factor of
$
k!\left( \left( {n\over k} \right) ! \right) ^k
$ that is greater than $1$
can be matched with a corresponding
equal or larger factor of
$\ell !\left( \left( {n\over \ell } \right) ! \right) ^\ell $.

For the second inequality,
because $\ell $ is the greatest nontrivial divisor of
$n$, two factors $2$ of
$\ell !\left( \left( {n\over \ell } \right) ! \right) ^\ell $
can be matched with a factor
greater than or equal to $4$ of
$(n-1)!$
and each of the remaining factors
greater than $1$ can be matched with
a corresponding
equal or larger factor of
$(n-1)!$.
\qed

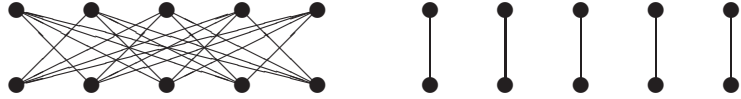
\begin{figure}

\centerline{
\unitlength 1mm 
\linethickness{0.4pt}
\ifx\plotpoint\undefined\newsavebox{\plotpoint}\fi 
\begin{picture}(101,16)(0,0)
\put(45,5){\circle*{2}}
\put(45,15){\circle*{2}}
\put(35,5){\circle*{2}}
\put(35,15){\circle*{2}}
\put(25,5){\circle*{2}}
\put(25,15){\circle*{2}}
\put(70,5){\circle*{2}}
\put(90,5){\circle*{2}}
\put(15,5){\circle*{2}}
\put(70,15){\circle*{2}}
\put(90,15){\circle*{2}}
\put(15,15){\circle*{2}}
\put(80,5){\circle*{2}}
\put(100,5){\circle*{2}}
\put(60,5){\circle*{2}}
\put(5,5){\circle*{2}}
\put(80,15){\circle*{2}}
\put(100,15){\circle*{2}}
\put(60,15){\circle*{2}}
\put(5,15){\circle*{2}}
\put(70,5){\line(0,1){10}}
\put(90,5){\line(0,1){10}}
\put(80,5){\line(0,1){10}}
\put(100,5){\line(0,1){10}}
\put(60,5){\line(0,1){10}}
\put(25,15){\line(-2,-1){20}}
\put(35,15){\line(-2,-1){20}}
\put(45,15){\line(-2,-1){20}}
\put(25,5){\line(-2,1){20}}
\put(35,5){\line(-2,1){20}}
\put(45,5){\line(-2,1){20}}
\put(5,5){\line(1,1){10}}
\put(15,5){\line(1,1){10}}
\put(25,5){\line(1,1){10}}
\put(35,5){\line(1,1){10}}
\put(15,5){\line(-1,1){10}}
\put(25,5){\line(-1,1){10}}
\put(35,5){\line(-1,1){10}}
\put(45,5){\line(-1,1){10}}
\put(45,15){\line(-4,-1){40}}
\put(45,5){\line(-4,1){40}}
\put(5,5){\line(3,1){30}}
\put(15,5){\line(3,1){30}}
\put(5,15){\line(3,-1){30}}
\put(15,15){\line(3,-1){30}}
\end{picture}
}

\caption{The standard example $S_5 $ and the disjoint union $5C_2 $ of
five $2$-chains.
}
\label{st_and_2C}

\end{figure}

\begin{lem}
\label{firstauffaecherung}

Let $P$ be an ordered set of height $1$
such that the following hold.
\begin{enumerate}
\item
$P$ has
$w\geq 3$ minimal elements.
\item
There is an $\ell \in \{ 1, \ldots , w-1\} $ such that
every maximal element is above exactly $\ell $
minimal elements.
\item
$P$ does not contain any nontrivial order-autonomous antichains.
\item
$\{ \Phi |_{R_0} :\Phi \in {\rm Aut} (P)\} $
contains $A_w (R_0 )$.
\end{enumerate}
Then $R_1 $ has exactly $\pmatrix{ w\cr \ell \cr } $ elements.
For $\ell \geq 2$,
$P$ is isomorphic to the ordered set that consists of the
singleton subsets and the $\ell $-element subsets of $R_0 $
ordered by inclusion.
For $\ell =1$, $P$ is isomorphic to the
pairwise disjoint union $wC_2 $
(see Figure \ref{st_and_2C})
of $w$ chains with $2$ elements each.

\end{lem}

{\bf Proof.}
Let $x\in R_1 $. Let $A:=\downarrow x\cap R_0 $
and let $B$ be any
$\ell $-element subset of
$R_0 $.

We first claim that
there is a
$\Phi _\beta \in {\rm Aut} (P)$ such
that
$\Phi _\beta [A]=B$.
First consider the case $\ell \geq 2$.
In this case, there are
$|A\setminus B|$
pairwise disjoint transpositions
whose composition $\gamma $ maps
$A\setminus B$ to $B\setminus A$ and
which leaves all other points
fixed.
Hence $\gamma $ maps $A$ to $B$.
If
$|A\setminus B|$ is even, then
$\gamma \in A_w (R_0 )$ and we set $\beta :=\gamma $.
If
$|A\setminus B|$ is odd, let
$\tau $ be a transposition that interchanges
two elements of $B$.
Then $\beta := \tau \circ \gamma $ maps $A$ to $B$
and
$\beta \in A_w (R_0 )$.
By assumption, there is a $\Phi _\beta \in {\rm Aut} (P)$
such that $\Phi _\beta |_{R_0 } =\beta $,
which proves the claim for $\ell \geq 2$.

For $\ell =1$, let $A=\{ a\} $, let $B=\{ b\} $,
let $\gamma :=(ab)$ and let
$\tau $ be a transposition that interchanges
two elements of $R_0 \setminus B$.
Then $\beta :=
\tau \circ \gamma $ is even and maps $A$ to $B$.
By assumption, there is a $\Phi _\beta \in {\rm Aut} (P)$
such that $\Phi _\beta |_{R_0 } =\beta $,
which proves the claim for $\ell =1$.

Because
$\Phi _\beta \in {\rm Aut} (P)$ and
$\Phi _\beta [A]=B$, we conclude that
$\downarrow \Phi _\beta (x)\cap R_0 =B$.
Because $B$ was arbitrary, for every
$\ell $-element subset $L$ of $R_0 $, there is an
$x_L \in R_1 $ such that
$\downarrow x_L \cap R_0 =L$.
Because $P$ does not contain any nontrivial
order-autonomous antichains,
this element $x_L $ is unique.
Because every element of $R_1 $ is above exactly $\ell $
elements of $R_0 $, there are no further elements in $R_1 $.

For $\ell =1$, the claimed isomorphism is clear.
For $\ell \geq 2$,
the claimed isomorphism maps
every $z\in R_0 $ to
$\{ z\} $ and every
$x\in R_1 $ to
$\downarrow x \cap R_0 $.
The claim about the number of
elements in $R_1 $ follows easily.
\qed

\begin{define}

(See Figure \ref{st_and_2C}.)
Let $w\geq 3$. We define the ordered set $S_w $ to be
the ordered set consisting of the $1$-element subsets
and the $w-1$-element subsets of $\{ 1, \ldots , w\} $
ordered by inclusion.

\end{define}

Recall that an ordered set $P$ is {\bf coconnected}
iff
there are no two nonempty subsets $A$ and $B$ such that
$A<B$ and $P=A\cup B$.

\begin{lem}
\label{maxlock2lev}

Let $P$ be a coconnected ordered set of height $1$ and width $w\geq 3$
such that
$\{ \Phi |_{R_0} :\Phi \in {\rm Aut} (P)\} $
contains $A_w (R_0 )$.
Then $P$ is
isomorphic to
an ordered set $S_w $ or an ordered set $wC_2 $.
In particular $|{\rm Aut} (P)|=w!$.

\end{lem}

{\bf Proof.}
Recall that
automorphisms map order-autonomous antichains to
order-autonomous antichains.
If
$R_0 $ were to contain a nontrivial order autonomous antichain, then,
by Lemma \ref{noverklem}, we conclude that, for $w\geq 6$, we have
$|\{ \Phi |_{R_0} :\Phi \in {\rm Aut} (P)\} |<(w-1)!<{1\over 2} w!.$
For $w=4$, the same argument leads to
$|\{ \Phi |_{R_0} :\Phi \in {\rm Aut} (P)\} |\leq 8<{1\over 2} 4!$.
Because
$\{ \Phi |_{R_0} :\Phi \in {\rm Aut} (P)\} $
contains $A_w (R_0 )$, which has ${1\over 2} w!$ elements,
we conclude that $R_0 $ cannot contain any order-autonomous
antichains.
Moreover, because
$\{ \Phi |_{R_0} :\Phi \in {\rm Aut} (P)\} $
contains $A_w (R_0 )$
and $P$ has width $w$, we conclude that $|R_0 |=w$.

Because $P$ is coconnected, there are an
$\ell \in \{ 1, \ldots , w-1\} $ and an
$x\in R_1 $ such that $|\downarrow x\cap R_0 |=\ell $.
Obtain $P^* $ from $P$ by selecting exactly one element
from every order-autonomous antichain that has
either $0$ or exactly $\ell $ strict lower bounds.
Then $P^* $
is coconnected, has
height $1$, width $w$, no nontrivial order-autonomous antichains,
$\{ \Phi |_{R_0} :\Phi \in {\rm Aut} (P^* )\} $
contains $A_w (R_0 )$,
and every non-minimal
element has exactly $\ell $ strict lower bounds.
By Lemma \ref{firstauffaecherung}, because $P^* $ has width $w$,
we obtain that $\ell \in \{ 1, w-1\} $
and that
$P$ is
isomorphic to
an ordered set $S_w $ or an ordered set $wC_2 $.
{\em A fortiori},
because $P$ has width $w$, we obtain that $P$ contains
no nontrivial order-autonomous antichains.
Hence $P^* =P$ which establishes the claimed isomorphism
as well as $|{\rm Aut} (P)|=w!$.
\qed

\begin{prop}
\label{maxlockdef}

Let $P$ be a coconnected ordered set width $w\geq 3$
such that there is a $k$ such that
$\{ \Phi |_{R_k} :\Phi \in {\rm Aut} (P)\} $
contains $A_w (R_k )$.
Then, for every $j$ such that $R_{j+1} \not= \emptyset $,
we have that $R_j\cup R_{j+1} $ is isomorphic to
$S_w $ or $wC_2 $.
In particular, $|{\rm Aut} (P)|=w!$.

\end{prop}

{\bf Proof.}
First consider the case that there are elements of rank greater than $k$.
Then
$|R_k |=w$ and, because $P$ is coconnected,
$R_k \cup R_{k+1} $ is a
coconnected ordered set of height $1$ and width $w$
such that
$\{ \Phi |_{R_0} :\Phi \in {\rm Aut} (P)\} $
contains $A_w (R_0 )$.
By Lemma \ref{maxlock2lev},
we have that $R_k\cup R_{k+1} $ is isomorphic to
$S_w $ or $wC_2 $.

By Lemma \ref{allbutoneac} applied to
$R_k\cup R_{k+1} $, we obtain that
every $\Phi |_{R_{k+1}} $ is uniquely determined by
$\Phi |_{R_{k}} $.
Because
$\{ \Phi |_{R_{k}} :\Phi \in {\rm Aut} (P)\} $
contains $A_w (R_{k} )$,
$\{ \Phi |_{R_{k+1}} :\Phi \in {\rm Aut} (P)\} $
contains $A_w (R_{k+1} )$.
Inductively, with $h$ denoting the height of $P$,
we infer that, for every
$j\in \{ k, \ldots , h-1\} $,
the ordered set $R_j\cup R_{j+1} $ is isomorphic to
$S_w $ or $wC_2 $.

Now, independent of whether
$k=h$ or $k<h$,
the
elements of rank $k$ are
the elements of dual rank $h-k$.
The dual of the preceding argument finishes the proof.
\qed

\begin{define}

An ordered set $P$ of width $w\geq 2$
such that, for every $j$ such that $R_{j+1} \not= \emptyset $,
we have that $R_j\cup R_{j+1} $ is isomorphic to
$S_w $ or $wC_2 $
will be called {\bf max-locked}.

\end{define}

\begin{remark}

{\rm
Note that max-locked ordered sets such that,
for every $j$ such that $R_{j+1} \not= \emptyset $,
we have that $R_j\cup R_{j+1} $ is isomorphic to
$wC_2 $ need not be unions of pairwise disjoint chains:
There can be $k\geq 0$ and $\ell \geq 2$ such that
$R_k <R_{k+\ell } $.
}

\end{remark}

\section{Interdependent Orbit Unions}
\label{orbclustsec}

The proof of Proposition \ref{maxlockdef}
shows how the action of the automorphism group on
a single set of elements of rank $k$ can, in
natural fashion, ``transmit vertically" though the whole ordered set.
In the situation of Proposition \ref{maxlockdef}, the sets $R_k $ happen to be
orbits
(see Definition \ref{grouporb} below) of $P$.
In general, the ``transmission" must
focus on the orbits, not the sets $R_k $, and it
can ``transmit" the action of the automorphism group
on one orbit $O$ to orbits that have no points that are
comparable to any element of $O$.
It should be noted that interdependent orbit unions
(see Definition \ref{interdeporbundef} below)
in graphs have also proven useful
in the set reconstruction of certain graphs, see \cite{SchrSetRec}.
Although the presentation up to Proposition \ref{getallfromiou}
translates directly to and from the corresponding results in
\cite{SchrSetRec}, all proofs are included to keep this presentation self-contained.

\begin{define}
\label{grouporb}

(Compare with Definition 8.1 in \cite{SchrSetRec}.)
Let $P$ be an ordered set, let $G$ be a subgroup of
${\rm Aut} (P)$ and let
$x\in P$.
Then the set
$G\cdot x :=\{ \Phi (x):\Phi \in G \} $
is called
the {\bf orbit of $x$ under the action of $G$} or the
{\bf $G$-orbit of $x$}.
Explicit mention of $G$ or $x$
can be dropped when there is only one group under consideration or when
specific knowledge of $x$ is not needed.
When no group $G$ is explicitly
mentioned at all, we assume by default that
$G={\rm Aut} (P)$.
The group generated by a single automorphism is denoted $\langle \Phi \rangle $.

\end{define}

\begin{figure}

\centerline{
\unitlength 1mm 
\linethickness{0.4pt}
\ifx\plotpoint\undefined\newsavebox{\plotpoint}\fi 
\begin{picture}(107,36.5)(0,0)
\put(5,5){\circle*{2}}
\put(5,25){\circle*{2}}
\put(105,5){\circle*{2}}
\put(105,25){\circle*{2}}
\put(5,15){\circle*{2}}
\put(5,35){\circle*{2}}
\put(105,15){\circle*{2}}
\put(105,35){\circle*{2}}
\put(15,5){\circle*{2}}
\put(15,25){\circle*{2}}
\put(95,5){\circle*{2}}
\put(95,25){\circle*{2}}
\put(15,15){\circle*{2}}
\put(15,35){\circle*{2}}
\put(95,15){\circle*{2}}
\put(95,35){\circle*{2}}
\put(25,5){\circle*{2}}
\put(25,25){\circle*{2}}
\put(85,5){\circle*{2}}
\put(85,25){\circle*{2}}
\put(25,15){\circle*{2}}
\put(25,35){\circle*{2}}
\put(85,15){\circle*{2}}
\put(85,35){\circle*{2}}
\put(45,5){\circle*{2}}
\put(65,5){\circle*{2}}
\put(55,5){\circle*{2}}
\put(55,5){\circle*{2}}
\put(65,5){\circle*{2}}
\put(45,5){\circle*{2}}
\put(5,5){\line(0,1){10}}
\put(5,15){\line(0,1){10}}
\put(85,15){\line(0,1){10}}
\put(5,25){\line(0,1){10}}
\put(105,5){\line(0,1){10}}
\put(105,25){\line(0,1){10}}
\put(15,5){\line(0,1){10}}
\put(15,15){\line(0,1){10}}
\put(95,15){\line(0,1){10}}
\put(15,25){\line(0,1){10}}
\put(95,5){\line(0,1){10}}
\put(95,25){\line(0,1){10}}
\put(25,5){\line(0,1){10}}
\put(25,15){\line(0,1){10}}
\put(105,15){\line(0,1){10}}
\put(25,25){\line(0,1){10}}
\put(85,5){\line(0,1){10}}
\put(85,25){\line(0,1){10}}
\put(25,5){\line(-2,1){20}}
\put(25,15){\line(-2,1){20}}
\put(105,15){\line(-2,1){20}}
\put(5,15){\line(2,1){20}}
\put(85,15){\line(2,1){20}}
\put(25,25){\line(-2,1){20}}
\put(85,5){\line(2,1){20}}
\put(85,25){\line(2,1){20}}
\put(45,5){\line(-4,1){40}}
\put(65,5){\line(4,1){40}}
\put(55,5){\line(-4,1){40}}
\put(55,5){\line(4,1){40}}
\put(65,5){\line(-4,1){40}}
\put(45,5){\line(4,1){40}}
\put(5,5){\line(1,1){10}}
\put(5,15){\line(1,1){10}}
\put(85,15){\line(1,1){10}}
\put(25,15){\line(-1,1){10}}
\put(105,15){\line(-1,1){10}}
\put(5,25){\line(1,1){10}}
\put(105,5){\line(-1,1){10}}
\put(105,25){\line(-1,1){10}}
\put(15,5){\line(1,1){10}}
\put(15,15){\line(1,1){10}}
\put(95,15){\line(1,1){10}}
\put(15,15){\line(-1,1){10}}
\put(95,15){\line(-1,1){10}}
\put(15,25){\line(1,1){10}}
\put(95,5){\line(-1,1){10}}
\put(95,25){\line(-1,1){10}}
\put(3,5){\makebox(0,0)[rc]{\footnotesize $A$}}
\put(3,25){\makebox(0,0)[rc]{\footnotesize $C$}}
\put(107,5){\makebox(0,0)[lc]{\footnotesize $\widetilde{A} $}}
\put(107,25){\makebox(0,0)[lc]{\footnotesize $\widetilde{C} $}}
\put(67,5){\makebox(0,0)[lt]{\footnotesize $M$}}
\put(3,15){\makebox(0,0)[rc]{\footnotesize $B$}}
\put(3,35){\makebox(0,0)[rc]{\footnotesize $D$}}
\put(107,15){\makebox(0,0)[lc]{\footnotesize $\widetilde{B} $}}
\put(107,35){\makebox(0,0)[lc]{\footnotesize $\widetilde{D} $}}
\put(45,5){\line(5,1){50}}
\put(55,5){\line(5,1){50}}
\put(65,5){\line(2,1){20}}
\multiput(5,25)(.2693602694,.0336700337){297}{\line(1,0){.2693602694}}
\multiput(105,25)(-.2693602694,.0336700337){297}{\line(-1,0){.2693602694}}
\multiput(15,25)(.2693602694,.0336700337){297}{\line(1,0){.2693602694}}
\multiput(95,25)(-.2693602694,.0336700337){297}{\line(-1,0){.2693602694}}
\multiput(25,25)(.2693602694,.0336700337){297}{\line(1,0){.2693602694}}
\multiput(85,25)(-.2693602694,.0336700337){297}{\line(-1,0){.2693602694}}
\put(15,5){\oval(23,3)[]}
\put(95,5){\oval(23,3)[]}
\put(55,5){\oval(23,3)[]}
\put(15,15){\oval(23,3)[]}
\put(95,15){\oval(23,3)[]}
\put(15,25){\oval(23,3)[]}
\put(95,25){\oval(23,3)[]}
\put(15,35){\oval(23,3)[]}
\put(95,35){\oval(23,3)[]}
\end{picture}
}

\caption{An ordered set
$P$
with
orbits marked with ovals.
}
\label{transmit_drive}

\end{figure}
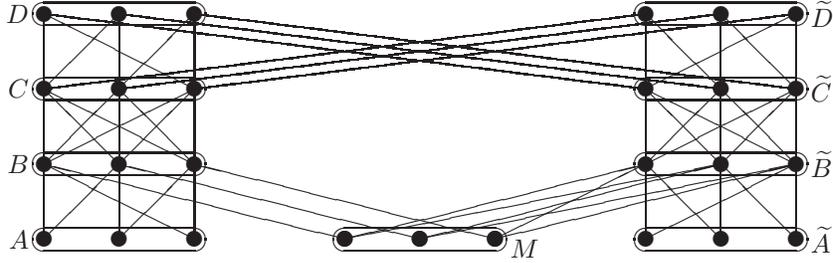

Note that,
if a strict subset $Q\subset P$ was obtained by
removing a union of ${\rm Aut} (P)$-orbits, then
${\rm Aut} (P)$-orbits
that are contained in $Q$
can be strictly contained in
${\rm Aut} (Q)$-orbits:
The ${\rm Aut} (P)$-orbits of the ordered set $P$ in
Figure \ref{transmit_drive} are marked by ovals.
We can see that, for $X\in \{ A,B,C,D\} $,
the ${\rm Aut} (P)$-orbits $X$ and $\widetilde{X} $
are strictly contained in
the ${\rm Aut} (P\setminus M)$-orbit $X\cup \widetilde{X} $.
For this reason (further elaborated later
in Remark \ref{orbstrictcont}),
dictated orbit structures/orbit frames will be
useful for the representation of ${\rm Aut} (P)$ in Proposition \ref{getallfromiou}
and they are vital for the remaining investigation.

\begin{define}
\label{dictatedef}

(Compare with Definition 8.2 in \cite{SchrSetRec}.)
Let $P$ be an ordered set and let ${\cal D}$ be a partition of $P$
into antichains. Then ${\rm Aut} _{\cal D} (P)$ is the set of automorphisms
$\Phi :P\to P$ such that, for every $\Phi $-orbit
$O:=\langle \Phi \rangle \cdot x$ of $\Phi $, there is a
$D\in {\cal D}$ such that $O\subseteq D$.
In this context, the partition ${\cal D}$ is called a
{\bf dictated orbit structure} (\cite{SchrSetRec}) or
an {\bf orbit frame}\footnote{Although this is a revision of
the author's own terminology, the investigation in this
paper and the fact that the elements of ${\cal D}$ need not be orbits
make it clear that the terminology ``orbit frame" is more appropriate.
The author thanks Frank a Campo for this suggestion.}
for $P$, and
the pair $(P,{\cal D})$ is called a {\bf structured ordered set}.
${\rm Aut} _{\cal D} (P)$-orbits will, more briefly, be called
{\bf ${\cal D}$-orbits}.
Sets $D\in {\cal D}$ will be called {\bf frames}.

The partition of $P$ into its ${\rm Aut} (P)$-orbits is called the
{\bf natural orbit structure/frame} of $P$, which will typically be denoted ${\cal N}$.
When working with the natural orbit structure/frame,
explicit indications of the automorphism set, usually via subscripts
or prefixes
${\cal D}$, will often be
omitted.

\end{define}

Clearly,
${\rm Aut} _{\cal D} (P)$ is a subgroup of the
automorphism group
${\rm Aut} (P)$.
Moreover,
${\rm Aut} _{\cal D} (P)\not= {\rm Aut} (P)$
iff there are an orbit $O$ of $P$ and a $D\in {\cal D}$
such that $O\cap D\not= \emptyset $ and
$O\not\subseteq D$.

The proof of
Proposition \ref{maxlockdef} has already shown
that, for an automorphism $\Phi $,
the values on a single orbit of $\Phi $
can completely determine
$\Phi $.
The relation of direct interdependence in Definition \ref{directdep}
below
provides a more detailed view of this observation as well as
of the
simple observation in
Lemma \ref{allbutoneac}.
Recall that $x\sim y$ denotes the fact that $x\leq y$ or $x\geq y$

\begin{define}
\label{directdep}

(Compare with Definition 8.3 in \cite{SchrSetRec}.)
Let
$(P,{\cal D})$ be a structured ordered set
and let $C,D$ be two ${\cal D}$-orbits
of $P$ such that there are a $c_1 \in C$
and a $d_1 \in D$ such that $c_1 <d_1 $.
We will write $C\upharpoonleft _{\cal D} D$
and $D\downharpoonright _{\cal D} C$
iff there are $c_2 \in C$ and a $d_2 \in D$
such that $c_2 \not\sim d_2 $.
In case $C\upharpoonleft _{\cal D} D$ or $C\downharpoonright _{\cal D} D$, we write
$C\updownharpoons _{\cal D} D$ and say that
$C$ and $D$ are
{\bf directly interdependent}.

\end{define}

Figure \ref{transmit_drive} shows
how
(connection through)
direct interdependence can
allow one orbit to determine the values of
automorphisms on many other orbits:
In the ordered set in
Figure \ref{transmit_drive}, we have
$A\upharpoonleft _{\cal N} B\downharpoonright _{\cal N} M
\upharpoonleft _{\cal N} \widetilde{B}
\downharpoonright _{\cal N} \widetilde{A}$ and the
values of any automorphism on
$A\cup B\cup M\cup \widetilde{B}
\cup \widetilde{A}$
are determined the by automorphism's values on $A$.
Note, however, that even direct interdependence
of two orbits $C\updownharpoons _{\cal N} D$
does not mean that
automorphisms are determined by the values on either
orbit
$C$ or $D$:
Consider a $6$-crown in which every element is replaced with an
order-autonomous $2$-antichain:
The minimal elements form an orbit, as do the maximal elements, the two
orbits are directly interdependent, but no automorphism is
completely determined solely
by its values on the minimal (or the maximal) elements.

We now turn to interdependent orbit unions, which
are the unions of the connected components of the orbit graph
defined in Definition \ref{orbitgraph} below.
Orbit graphs themselves will take center stage starting in Section
\ref{orbgraphsec}.

\begin{define}
\label{orbitgraph}

Let
$(P, {\cal D} )$
be
a structured ordered set.
The {\bf orbit graph}
${\cal O} (P, {\cal D})$
of
$P$
is defined to be the graph whose vertex set is
the set of all ${\cal D}$-orbits
such that two ${\cal D}$-orbits
$D_1 $ and $D_2 $ are adjacent iff
$D_1 \updownharpoons _{\cal D} D_2 $.

\end{define}

Although our main focus will
ultimately be on orbit graphs
for structured ordered sets in which every $D\in {\cal D}$ is a ${\cal D} $-orbit,
also see part \ref{structset3}
of Definition \ref{structset} below,
until that time, we must keep in mind that the
$D\in {\cal D}$ are frames which need not be orbits themselves.

\begin{define}
\label{interdeporbundef}

(Compare with Definition 8.5 in \cite{SchrSetRec}.)
Let
$(P, {\cal D} )$
be
a structured ordered set and let ${\cal O} (P,{\cal D} )$ be its orbit graph.
If ${\cal E} $ is a connected component of
${\cal O} (P,{\cal D} )$,
then we will call the set $\bigcup {\cal E}\subseteq P$ an
{\bf interdependent ${\cal D}$-orbit union}.

\end{define}

\begin{define}

(Compare with Definition 8.6 in \cite{SchrSetRec}.)
Let
$(P,{\cal D})$ be a structured ordered set
and let
$Q\subseteq P$ such that, for all
${\cal D}$-orbits $D$,
we have $D\subseteq Q $ or $D\cap Q =\emptyset $.
The
orbit frame for $Q$
{\bf induced by ${\cal D}$}, denoted ${\cal D} |Q $,
is defined to be the set of all ${\cal D}$-orbits
that are contained in $Q $.

For the natural orbit
frame ${\cal N}$ for $P$,
the
orbit frame for $Q$ induced by ${\cal N} $
will be called the {\bf naturally required}
orbit
frame
${\cal N}|Q$.

\end{define}

Note that, for induced orbit
frames
${\cal D} |Q$, every
$D\in {\cal D} |Q$ is an orbit.
Proposition \ref{unionplacement}
below
now shows how
interdependent ${\cal D}$-orbit unions
reside in an ordered set.
It also lays the groundwork for representing
automorphisms through certain automorphisms on the
non-singleton interdependent orbit unions in
Proposition \ref{getallfromiou}.

\begin{prop}
\label{unionplacement}

(Compare with Proposition 8.7 in \cite{SchrSetRec}.)
Let
$(P,{\cal D})$ be a structured ordered set
and let
$U$ be an interdependent ${\cal D}$-orbit union.
Then the following hold.
\begin{enumerate}
\item
\label{unionplacement0}
For all $x\in P\setminus U$ and all
$C\in {\cal D}|U$, the following hold.
\begin{enumerate}
\item
\label{unionplacement1}
If there is a $c\in C$ such that $c<x$, then $C<x$.
\item
\label{unionplacement2}
If there is a $c\in C$ such that $c>x$, then $C>x$.
\end{enumerate}
\item
\label{unionplacement3a}
For every $\Phi \in {\rm Aut} _{\cal D} (P)$, we have that
$\Phi |_U \in {\rm Aut} _{{\cal D}|U} (U)$.
In particular, this means that the
${\cal D}|U$-orbits
are just the sets in
${\cal D}|U$.
\item
\label{unionplacement3b}
For
every $\Psi \in {\rm Aut} _{{\cal D}|U} (U)$, the function
$\Psi ^P (x) := \cases{
\Psi (x); & if $x\in U$, \cr
x; & if $x\in P\setminus U$,
} $
is an automorphism of $P$.
\end{enumerate}

\end{prop}

{\bf Proof.}
To prove part \ref{unionplacement1},
let $x\in P\setminus U$ and
$C\in {\cal D}|U$ be so that
there is a $c\in C$ such that $c<x$.
Let $X$ be the ${\cal D}$-orbit of $x$.
Because
$x\not\in U$, we must have
$X\not\updownharpoons _{\cal D} C$.
Because $C\ni c<x\in X$, we must have that
$C<X$ and hence $C<x$.

Part \ref{unionplacement2} is proved dually.

Part \ref{unionplacement3a}
follows directly from the definitions.

To prove part \ref{unionplacement3b},
let $\Psi \in {\rm Aut} _{{\cal D}|U} (U)$.
Clearly, $\Psi ^P $ is bijective.
To prove that $\Psi ^P $ is
order-preserving, let $x<y$.
If $x,y$ are both in $U$ or both in $P\setminus U$, we obtain
$\Psi ^P (x)<\Psi ^P (y)$.
In case $x\in P\setminus U$ and $y\in U$,
let $Y\in {\cal D}|U$ be so that $y\in Y$.
Then $\Psi (y)\in Y$.
By part \ref{unionplacement2},
we have that $x<Y$ and hence
$\Psi ^P (x)=x<Y\ni \Psi ^P (y)$.
The case in which $y\in P\setminus U$ and $x\in U$
is handled dually.
\qed

\begin{define}

(Compare with Definition 8.8 in \cite{SchrSetRec}.)
Let
$(P,{\cal D})$ be a structured ordered set
and let
$U$ be an interdependent ${\cal D}$-orbit union.
We define ${\rm Aut} _{{\cal D}|U} ^P (U)$ to be the set of
automorphisms $\Psi ^P \in {\rm Aut} (P)$ as in
part \ref{unionplacement3b} of Proposition \ref{unionplacement}.

\end{define}

Let $P$ be an ordered set, let ${\cal D}$ be
an orbit frame
for $P$ and let
$U, U'$ be disjoint interdependent ${\cal D}$-orbit unions.
Then, clearly, for
$\Psi ^P \in {\rm Aut} _{{\cal D}|U} ^P (U)$
and
$\Phi ^P \in {\rm Aut} _{{\cal D}|U'} ^P (U')$, we have
$\Psi ^P \circ \Phi ^P =\Phi ^P \circ \Psi ^P $.
Hence the following definition is sensible.

\begin{define}

(Compare with Definition 8.9 in \cite{SchrSetRec}.)
Let $P$ be an ordered set and let
${\cal A}_1 , \ldots , {\cal A} _z \subseteq {\rm Aut} (P)$ be
sets of automorphisms such that, for all pairs of distinct
$i,j\in \{ 1, \ldots , z\} $,
all $\Phi _i \in {\cal A} _i $ and all $\Phi _j \in {\cal A} _j $,
we have $\Phi _i \circ \Phi _j = \Phi _j \circ \Phi _i $.
We define $\bigcirc _{j=1} ^z
{\cal A} _j $
to be the set of
compositions $\Psi _1 \circ \cdots \circ \Psi _z $
such that, for $j=1, \ldots , z$, we have
$\Psi _j \in
{\cal A}_j $.

\end{define}

\begin{prop}
\label{getallfromiou}

(Compare with Proposition 8.10 in \cite{SchrSetRec}.)
Let $P$ be an
ordered set
with natural orbit
frame ${\cal N}$
and let $U_1 , \ldots , U_z $
be the non-singleton interdependent orbit unions of $P$.
Then
${\rm Aut} (P)=\bigcirc _{j=1} ^z
{\rm Aut} _{{\cal N} |U_j } (U_j )$, and consequently
$
|{\rm Aut} (P)|=\prod _{j=1} ^z \left|{\rm Aut} _{{\cal N} |U_j } ^P (U_j )\right|
$.

\end{prop}

{\bf Proof.}
The containment ${\rm Aut} (P)\supseteq \bigcirc _{j=1} ^z
{\rm Aut} _{{\cal N} |U_j } ^P (U_j )$
follows from part \ref{unionplacement3b} of Lemma \ref{unionplacement}.

By part \ref{unionplacement3a} of Lemma \ref{unionplacement},
for every $\Phi \in {\rm Aut} (P)$ and every $j\in \{ 1, \ldots , z\} $,
we have that
$\Phi |_{U_j } \in {\rm Aut} _{{\cal N}|U_j } (U_j)$.
Because $\Phi $ fixes
all points in $P\setminus \bigcup _{j=1} ^z U_j $,
we have
$\Phi =\Phi |_{U_1 } ^P \circ \cdots \circ \Phi |_{U_z } ^P $.
Hence
${\rm Aut} (P)\subseteq \bigcirc _{j=1} ^z
{\rm Aut} _{{\cal N} |U_j } ^P (U_j )$.
\qed

\begin{remark}
\label{orbstrictcont}

(Compare with Remark 8.11 in \cite{SchrSetRec}.)
{\rm
Although the
naturally required
orbit frame
${\cal N}|U$
may look more technical than natural, it is indispensable for
the representation in Proposition \ref{getallfromiou}.
Consider the ordered set in Figure \ref{transmit_drive}.
The natural interdependent
orbit unions in this ordered set are
$U_1 :=
A\cup B\cup M\cup \widetilde{B}\cup \widetilde{A} $
and
$U_2 :=
C\cup D\cup \widetilde{D}\cup \widetilde{C} $.
However, when considering $U_2 $
as an ordered set by itself, the
${\rm Aut} (U_2 )$-orbits
are $C\cup \widetilde{C}$ and
$D\cup \widetilde{D}$, whereas the
${\rm Aut} (P)$-orbits in $U_2 $
are
$C, D, \widetilde{C}$ and
$\widetilde{D}$.
Hence, we cannot use
the automorphism groups
${\rm Aut} (U_j )$
in place of their subgroups
${\rm Aut} _{{\cal N} |U_j } (U_j )$
in Proposition \ref{getallfromiou}.

The same effect was observed
in the introduction in the case that
$M$ is removed from the ordered set.
This is
why
orbit frames are also
crucial for the proof of
Theorem \ref{permgrtoordaut}.
\qex
}

\end{remark}

When bounding the number of automorphisms,
Proposition \ref{getallfromiou} allows us to focus our
efforts on ordered sets that
consist of a single interdependent orbit union $(U,{\cal D})$.
Even when not stated explicitly later,
any upper bound $|{\rm Aut} _{\cal D} (U)|\leq 2^{c|U|} $
with $c\leq 0.66$ (or, in light of
Theorem \ref{moreendos}, $c\leq 0.79$)
gives another class of ordered sets for which the
Automorphism Conjecture is confirmed.

We conclude this section by introducing ideas and terminology that will be
fundamental for the remainder of this paper.

\begin{define}

We will say that the structured ordered set
$(P,{\cal D} _P )$ is
{\bf (dually) isomorphic}
to the structured ordered set
$(Q,{\cal D} _Q )$ iff there is a (dual) isomorphism
$\Phi :P\to Q$ such that, for all
$D\in {\cal D} _P $, we have
$\Phi [D]\in {\cal D} _Q $.

\end{define}

\begin{define}
\label{structset}

Let $(P, {\cal D} )$ be a structured ordered set.
\begin{enumerate}
\item
\label{structset1}
We call $(P,{\cal D})$ an {\bf interdependent orbit union}
iff
$\bigcup {\cal D} $
is an interdependent ${\cal D}$-orbit union.

\item
\label{structset2}
We say $(P,{\cal D})$ is
{\bf without slack}
iff
none of the sets in ${\cal D}$ contains a nontrivial order-autonomous
antichain.

\item
\label{structset3}
We call $(P,{\cal D})$ {\bf tight}
iff
it is without slack and
each $D \in {\cal D} $ is a
${\cal D}$-orbit, that is,
it is without slack and
${\rm Aut} _{\cal D} (P)$ acts transitively on
every $D \in {\cal D} $.

\end{enumerate}

\end{define}

\begin{define}

Let $(P,{\cal D} )$ be a structured ordered set and let
$S\subseteq P$.
We define
$\Lambda _{\cal D} (S):=\{ \Phi |_S :\Phi \in {\rm Aut} _{\cal D} (P)\} $.

\end{define}

Our main focus will be on
interdependent orbit unions $(U,{\cal D})$ such that
no $\Lambda _{\cal D} (D)$ contains $A_{|D|} (D)$.
Because we will first focus on indecomposable ordered sets, we
will assume that $(U,{\cal D})$ is without slack.
Because we can
always refine the
orbit frame, we
are free to assume that $(U,{\cal D})$ is tight.

\section{The Orbit Graph of an Interdependent Orbit Union}
\label{orbgraphsec}

To bound the number of
automorphisms on a tight interdependent orbit
union $(U, {\cal D} )$, we will ultimately perform
an induction on the number of
frames.
Because $(U, {\cal D} )$ is tight,
every $D\in {\cal D}$ is an orbit.
Hence the vertex set of
${\cal O} (U, {\cal D})$ is ${\cal D}$.
When a
frame $D_n \in {\cal D}$ is removed,
the orbit graph can become disconnected and
the resulting structured ordered set may contain
nontrivial order-autonomous antichains.
We start with some notation that
assures that the orbits with indices smaller than $n$
form a component of the resulting
graph
${\cal O} (U, {\cal D})-D_{n} $, and
which allows easy reference to orbits with
certain properties.
This notation
will be used throughout the remaining sections.

\begin{nota}
\label{standardOGnotation}

Throughout,
$(U, {\cal D} )$
will be
a
tight
interdependent orbit union with $|{\cal D} |\geq 3$,
with the labeling of the elements of
${\cal D} =\{ D_1 , \ldots , D_m \} $
and $n\geq 3$ chosen
so that the following hold.
\begin{enumerate}
\item
$\{ D_1 , \ldots , D_{n-1} \} $
is a connected component of
the graph ${\cal O} (U, {\cal D})-D_{n} $.
\item
The orbits that are directly interdependent with
$D_n $ are $D_s , \ldots , D_{n-1} $ and $D_r , \ldots , D_m $.
\item
The orbits that contain nontrivial
$U\setminus \bigcup _{i=n} ^m D_i $-order-autonomous antichains are
$D_t , \ldots , D_{n-1} $.
\end{enumerate}

\end{nota}

Note that
Notation \ref{standardOGnotation} allows for $n=m$, that is, for
$D_n $ to not be a cutvertex of
${\cal O} (U, {\cal D})$.
Moreover, because
$(U, {\cal D} )$
is an
interdependent orbit union,
${\cal O} (U, {\cal D})$ is connected, so $s\leq n-1$.
For the other parameters,
$r> m$ or $t>n-1$ shall indicate that there are no orbits
as described via $r$ or $t$.

To easily refer to the (possibly trivial)
$U\setminus \bigcup _{i=n} ^m D_i $-order-autonomous
antichains in the orbits
$D_s , \ldots , D_{n-1} $,
we introduce the following notation.

\begin{nota}
\label{blocksinDj}

For every $j\in \{ s, \ldots , n-1\} $,
we let $A_1 ^j , \ldots , A_{\ell _j} ^j $ be the
maximal
$U\setminus \bigcup _{i=n} ^m D_i $-order-autonomous antichains
that partition $D_j $. We set
${\cal A} ^j :=\{ A_1 ^j , \ldots , A_{\ell _j} ^j \} $.
Moreover,
for every $i\in \{ 1, \ldots , \ell _j \} $, we choose
a fixed element $a_i ^j \in A_i ^j $.
Note that, for $s\leq j\leq t-1$, the sets $A_i ^j $ are singletons.

\end{nota}

The fact that automorphisms map maximal order-autonomous antichains to
maximal order-autonomous antichains motivates the definition below, which
will frequently be used and which leads to our first insights.

\begin{define}

Let $S$ be a set, let ${\cal A} =\{ A_1 , \ldots , A_\ell \} $
be a partition of $S$ and let $\Phi :S\to S$ be a permutation of $S$.
We say that {\bf $\Phi $ respects the partition ${\cal A}$}
iff, for every $i\in \{ 1, \ldots , \ell \} $, there is
a $j\in \{ 1, \ldots , \ell \} $ such that $\Phi [A_i ]=A_j $.

\end{define}

\begin{lem}
\label{pruneorbit1}

For every $j\in \{ s, \ldots n-1\} $ we have that
$\ell _j >1$ and
every
$\Phi \in {\rm Aut} _{\cal D} (U)$
respects the partition
${\cal A} ^j $
of $D_j $.
Moreover,
for all $i,k\in \{ 1, \ldots , \ell _j \} $,
we have
$\left| A_i ^j \right| =\left| A_k ^j \right| $.

\end{lem}

{\bf Proof.}
Let
$j\in \{ s, \ldots , n-1\} $.
Because $n\geq 3$
and
$\{ D_1 , \ldots , D_{n-1} \} $
is a connected component of
the graph ${\cal O} (U, {\cal D})-D_{n} $, we have that
$D_j $ is directly interdependent with another
$D_{j'} $ with $j'\in \{ 1, \ldots , n-1\} \setminus \{ j\} $.
Consequently,
$D_j $ itself is not a
$U\setminus \bigcup _{i=n} ^m D_i $-order-autonomous antichain.
Hence $\ell _j >1$.

Let
$\Phi \in {\rm Aut} _{\cal D} (U)$.
Because
$\Phi |_{U\setminus \bigcup _{i=n} ^m D_i } \in
{\rm Aut} _{{\cal D} \setminus \{ D_n , \ldots , D_m \} }
(U\setminus \bigcup _{i=n} ^m D_i )$, we have
$\Phi [D_j ]=D_j $,
and because automorphisms map
maximal order-autonomous antichains to
maximal order-autonomous antichains,
$\Phi $ must map every $A_i ^j $ to another
$A_k ^j $, that is,
$\Phi $
respects
${\cal A} ^j $.

Because $(U, {\cal D})$ is tight,
for any $i,k\in \{ 1, \ldots , \ell _j \} $,
there is a
$\Phi \in {\rm Aut} _{\cal D} (U)$ such that
$\Phi \left( a_i ^j \right) =a_k ^j $.
Consequently, because
$\Phi $
respects
${\cal A} ^j $, we have
$\Phi \left[ A_i ^j \right] =A_k ^j $, and hence
$\left| A_i ^j \right| =\left| A_k ^j \right| $.
\qed

\begin{lem}
\label{pruneorbit3}

For all
$j\in \{ t, \ldots , n-1\} $,
$i\in \{ 1, \ldots , \ell _j \} $
and distinct $x,y\in A_i ^j $, there
is a $d\in D_n $ such that $d$ is comparable to one of
$x$ and $y$, but not the other.

\end{lem}

{\bf Proof.}
Let
$j\in \{ t, \ldots , n-1\} $ and
$i\in \{ 1, \ldots , \ell _j \} $.
Because $D_j $ is not directly interdependent with any
$D_k $ with $k>n$, we have that
$A_i ^j $ is order-autonomous in $U\setminus D_n $.
Because, for any distinct $x,y\in A_i ^j $, the set
$\{ x,y\} $ is not order-autonomous in $U$,
there
must be a $d\in D_n $ such that $d$ is comparable to one of
$x$ and $y$, but not the other.
\qed

\vspace{.1in}

In our analysis, the order-autonomous antichains
$A_i ^j $ will be collapsed into
singletons. Hence we introduce the following and again establish some natural
properties.

\begin{define}
\label{phindef}

We define the
{\bf $D_{n} $-pruned and compacted}
ordered set $U_n :=\bigcup _{i=1} ^{s-1} D_i \cup \bigcup _{j=s} ^{n-1}
\left\{ a_i ^j : i\in \{ 1, \ldots , \ell _j \} \right\} $
and we define
\begin{eqnarray*}
{\cal D}_n
& := &
\{ D_j \cap U_n :j \in \{ 1, \ldots , n-1\} \}
\\
& = &
\{ D_j : j\in \{ 1, \ldots , s-1\} \} \cup
\left\{ \left\{ a_1 ^j , \ldots , a_{\ell _j } ^j \right\} :  j\in \{ s, \ldots , n-1\} \right\} .
\end{eqnarray*}
For every
$\Phi \in {\rm Aut} _{\cal D} (U)$, we define
the function $\Phi _n :U_n \to U_n  $
by
$\Phi _n |_{\bigcup _{j=1} ^{s-1} D_j } :=\Phi |_{\bigcup _{j=1} ^{s-1} D_j } $,
and by, for any $j\in \{ s, \ldots , n-1\} $ and
$i\in \{ 1, \ldots , \ell _j \} $, setting $\Phi _n \left( a_i ^j \right)
$ to be the unique element of $\Phi \left[ A_i ^j \right] \cap U_n $.

\end{define}

\begin{lem}
\label{pruneorbit2}

For every $\Phi \in {\rm Aut} _{\cal D} (U)$, we have that
$\Phi _n \in {\rm Aut } _{{\cal D}_n} (U_n )$.
Moreover,
$(U_n , {\cal D}_n )$ is a tight interdependent
orbit union.

\end{lem}

{\bf Proof.}
Clearly,
${\cal D}_n$ is an
orbit frame
for $U_n $.

Because $\{ D_1 , \ldots , D_{n-1} \} $
is a connected component of
the graph ${\cal O} (U, {\cal D})-D_{n} $,
$(U_n , {\cal D}_n )$ is an interdependent
orbit union.

Because $(U,{\cal D})$ is without slack and because we
choose exactly one element from each maximal
$U\setminus \bigcup _{i=n} ^m D_i $-order-autonomous antichain
to be in $U_n $, we obtain that $(U_n ,{\cal D}_n )$ is
an interdependent orbit union without slack.

For any
$\Phi \in {\rm Aut} _{\cal D} (U)$,
it follows from the definitions that
$\Phi _n \in {\rm Aut } _{{\cal D}_n} (U_n )$.
Because, for any $j\in \{ s, \ldots , n-1\} $ and
$x,y\in \{ 1, \ldots , \ell _j \} $, there is a
$\Phi \in {\rm Aut} _{\cal D} (U)$ with
$\Phi \left[ A_x ^j \right] =A_y ^j $, we conclude that
$(U_n , {\cal D}_n )$ is tight.
\qed

\vspace{.1in}

With the ``early orbits"
$D_1 , \ldots , D_{n-1} $ thus
discussed, we now turn to the remaining
orbits
$D_n , \ldots , D_m $ as well as the connection
between $D_n $ and the sets $A_i ^j $.

\begin{define}
\label{setQdef}

Let $Q:=\bigcup _{j=s} ^{m} D_j $,
let ${\cal E}_Q
:=
\{ A_i ^j : j=s, \ldots , n-1; i=1, \ldots , \ell _j \}
\cup \{ D_j :j=n, \ldots ,m \} $
and
let ${\cal D}_Q $ be the set of all ${\cal E}_Q $-orbits.

\end{define}

\begin{lem}
\label{Qonlystruct}

$(Q , {\cal D}_Q )$ is a tight
structured ordered set.

\end{lem}

{\bf Proof.}
First note that,
for $j\geq n$, all $D_i $ that are directly interdependent with
$D_j $ are contained in $\bigcup _{k=s} ^m D_k = Q$.
Therefore,
any $Q$-order-autonomous antichain in a set $D_j $ with $j\geq n$
would be $U$-order-autonomous.
Hence
for $j\geq n$, no $D_j $ contains a nontrivial
$Q$-order-autonomous antichain.

By
Lemma \ref{pruneorbit3},
for $j\in \{ s, \ldots , n-1\} $ and
$i\in \{ 1, \ldots , \ell _j \} $,
no $A_i ^j $ contains a nontrivial
$Q$-order-autonomous antichain.
We conclude that $(Q,{\cal E}_Q )$ is
a structured ordered set
without slack.
Because
${\cal D}_Q $ is the set of all ${\cal E}_Q $-orbits, we conclude that
$(Q,{\cal D}_Q )$ is a
tight
structured ordered set.
\qed

\begin{figure}

\centerline{
\unitlength 1mm 
\linethickness{0.4pt}
\ifx\plotpoint\undefined\newsavebox{\plotpoint}\fi 
\begin{picture}(113,44)(0,0)
\put(5,5){\line(0,1){10}}
\put(5,30){\line(0,1){10}}
\put(30,5){\line(0,1){10}}
\put(30,30){\line(0,1){10}}
\put(15,5){\line(0,1){10}}
\put(15,30){\line(0,1){10}}
\put(40,5){\line(0,1){10}}
\put(40,30){\line(0,1){10}}
\put(55,5){\line(0,1){10}}
\put(55,30){\line(0,1){10}}
\put(65,5){\line(0,1){10}}
\put(65,30){\line(0,1){10}}
\put(75,5){\line(0,1){10}}
\put(75,30){\line(0,1){10}}
\put(90,5){\line(0,1){10}}
\put(90,30){\line(0,1){10}}
\put(100,5){\line(0,1){10}}
\put(100,30){\line(0,1){10}}
\put(110,5){\line(0,1){10}}
\put(110,30){\line(0,1){10}}
\put(5,5){\circle*{2}}
\put(5,30){\circle*{2}}
\put(30,5){\circle*{2}}
\put(30,30){\circle*{2}}
\put(15,5){\circle*{2}}
\put(15,30){\circle*{2}}
\put(40,5){\circle*{2}}
\put(40,30){\circle*{2}}
\put(55,5){\circle*{2}}
\put(55,30){\circle*{2}}
\put(65,5){\circle*{2}}
\put(65,30){\circle*{2}}
\put(75,5){\circle*{2}}
\put(75,30){\circle*{2}}
\put(90,5){\circle*{2}}
\put(90,30){\circle*{2}}
\put(100,5){\circle*{2}}
\put(100,30){\circle*{2}}
\put(110,5){\circle*{2}}
\put(110,30){\circle*{2}}
\put(5,15){\circle*{2}}
\put(5,40){\circle*{2}}
\put(30,15){\circle*{2}}
\put(30,40){\circle*{2}}
\put(15,15){\circle*{2}}
\put(15,40){\circle*{2}}
\put(40,15){\circle*{2}}
\put(40,40){\circle*{2}}
\put(55,15){\circle*{2}}
\put(55,40){\circle*{2}}
\put(65,15){\circle*{2}}
\put(65,40){\circle*{2}}
\put(75,15){\circle*{2}}
\put(75,40){\circle*{2}}
\put(90,15){\circle*{2}}
\put(90,40){\circle*{2}}
\put(100,15){\circle*{2}}
\put(100,40){\circle*{2}}
\put(110,15){\circle*{2}}
\put(110,40){\circle*{2}}
\put(10,5){\oval(14,4)[]}
\put(10,30){\oval(14,4)[]}
\put(35,5){\oval(14,4)[]}
\put(35,30){\oval(14,4)[]}
\put(10,15){\oval(14,4)[]}
\put(10,40){\oval(14,4)[]}
\put(35,15){\oval(14,4)[]}
\put(35,40){\oval(14,4)[]}
\put(5,30){\line(5,2){25}}
\put(40,30){\line(-5,2){25}}
\put(30,40){\line(-3,-2){15}}
\put(15,40){\line(3,-2){15}}
\put(15,30){\line(5,2){25}}
\put(30,30){\line(-5,2){25}}
\multiput(40,40)(-.1178451178,-.0336700337){297}{\line(-1,0){.1178451178}}
\multiput(5,40)(.1178451178,-.0336700337){297}{\line(1,0){.1178451178}}
\put(65,5){\oval(24,4)[]}
\put(65,30){\oval(24,4)[]}
\put(100,5){\oval(24,4)[]}
\put(100,30){\oval(24,4)[]}
\put(65,15){\oval(24,4)[]}
\put(65,40){\oval(24,4)[]}
\put(100,15){\oval(24,4)[]}
\put(100,40){\oval(24,4)[]}
\multiput(55,30)(.1178451178,.0336700337){297}{\line(1,0){.1178451178}}
\multiput(110,30)(-.1178451178,.0336700337){297}{\line(-1,0){.1178451178}}
\put(90,40){\line(-5,-2){25}}
\put(75,40){\line(5,-2){25}}
\multiput(65,30)(.1178451178,.0336700337){297}{\line(1,0){.1178451178}}
\multiput(100,30)(-.1178451178,.0336700337){297}{\line(-1,0){.1178451178}}
\put(100,40){\line(-5,-2){25}}
\put(65,40){\line(5,-2){25}}
\multiput(75,30)(.1178451178,.0336700337){297}{\line(1,0){.1178451178}}
\multiput(90,30)(-.1178451178,.0336700337){297}{\line(-1,0){.1178451178}}
\multiput(110,40)(-.1515151515,-.0336700337){297}{\line(-1,0){.1515151515}}
\multiput(55,40)(.1515151515,-.0336700337){297}{\line(1,0){.1515151515}}
\multiput(100,40)(-.1515151515,-.0336700337){297}{\line(-1,0){.1515151515}}
\multiput(65,40)(.1515151515,-.0336700337){297}{\line(1,0){.1515151515}}
\multiput(55,30)(.1851851852,.0336700337){297}{\line(1,0){.1851851852}}
\multiput(110,30)(-.1851851852,.0336700337){297}{\line(-1,0){.1851851852}}
\put(75,30){\line(3,2){15}}
\put(90,30){\line(-3,2){15}}
\put(2,12){\dashbox{1}(41,6)[cc]{}}
\put(2,37){\dashbox{1}(41,6)[cc]{}}
\put(52,12){\dashbox{1}(61,6)[cc]{}}
\put(52,37){\dashbox{1}(61,6)[cc]{}}
\put(22,19){\makebox(0,0)[cb]{\footnotesize $D_n $}}
\put(22,44){\makebox(0,0)[cb]{\footnotesize $D_n $}}
\put(82,19){\makebox(0,0)[cb]{\footnotesize $D_n $}}
\put(82,44){\makebox(0,0)[cb]{\footnotesize $D_n $}}
\put(10,2){\makebox(0,0)[ct]{\footnotesize $A_1 ^{n-1} $}}
\put(10,27){\makebox(0,0)[ct]{\footnotesize $A_1 ^{n-1} $}}
\put(65,2){\makebox(0,0)[ct]{\footnotesize $A_1 ^{n-1} $}}
\put(65,27){\makebox(0,0)[ct]{\footnotesize $A_1 ^{n-1} $}}
\put(35,2){\makebox(0,0)[ct]{\footnotesize $A_2 ^{n-1} $}}
\put(35,27){\makebox(0,0)[ct]{\footnotesize $A_2 ^{n-1} $}}
\put(100,2){\makebox(0,0)[ct]{\footnotesize $A_2 ^{n-1} $}}
\put(100,27){\makebox(0,0)[ct]{\footnotesize $A_2 ^{n-1} $}}
\put(75,15){\line(-2,-1){20}}
\put(75,40){\line(-2,-1){20}}
\put(90,15){\line(2,-1){20}}
\put(90,40){\line(2,-1){20}}
\put(55,15){\line(1,-1){10}}
\put(55,40){\line(1,-1){10}}
\put(110,15){\line(-1,-1){10}}
\put(110,40){\line(-1,-1){10}}
\put(65,15){\line(1,-1){10}}
\put(65,40){\line(1,-1){10}}
\put(100,15){\line(-1,-1){10}}
\put(100,40){\line(-1,-1){10}}
\put(2,5){\makebox(0,0)[rc]{$S_4 ^\ddag $}}
\put(2,30){\makebox(0,0)[rc]{$S_4 ^\dagger $}}
\put(52,5){\makebox(0,0)[rc]{$S_6 ^\ddag $}}
\put(52,30){\makebox(0,0)[rc]{$S_6 ^\dagger $}}
\end{picture}
}

\caption{
Structured ordered sets induced on $D_{n-1} \cup D_n $ with $n=m$
that are not interdependent orbit unions.
The orbits are marked with ovals,
$D_n $ is marked with a dashed box,
and the partition $A_1 ^{n-1} \cup A_2 ^{n-1} $
is labeled.
}
\label{inconvenient}

\end{figure}
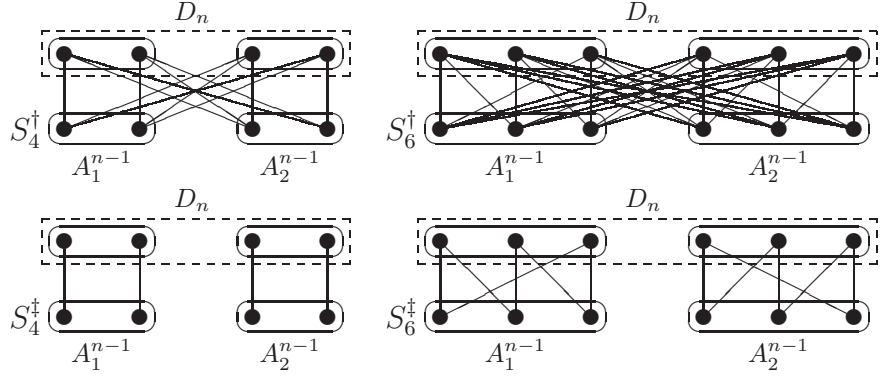

\vspace{.1in}

Lemma \ref{Qonlystruct} may feel a little unsatisfying compared to
Lemma \ref{pruneorbit2} in that $(Q,{\cal D} _Q )$
need not be an interdependent orbit union.
However, Figure \ref{inconvenient} gives examples that this really need not be the case:
Any of the
structured ordered sets given there could be a
structured ordered set $(Q,{\cal D} _Q )$ in the case
in which $m=n$,
$D_n $ is a pendant vertex, and $|U_n \cap D_{n-1} |=2$

\begin{define}

With
$\Phi _n $ as in Definition \ref{phindef},
we define
${\rm Aut} _{{\cal D} _Q } ^U (Q):=
\{ \Psi \in {\rm Aut} _{\cal D} (U): \Psi _n ={\rm id} _{U_n } \} $.
For every $\Delta \in
{\rm Aut } _{{\cal D}_Q } (Q)$, we define
$\Delta ^U $
by
$$\Delta ^U (x):=
\cases{
\Delta (x); & if $x\in Q$, \cr
x; & if $x\in U\setminus Q$.} $$

\end{define}

\begin{lem}
\label{AutDQUrepres}

${\rm Aut} _{{\cal D} _Q } ^U (Q)=
\{ \Psi \in {\rm Aut} _{\cal D} (U): \Psi _n ={\rm id} _{U_n } \}
=
\left\{ \Delta ^U :
\Delta \in {\rm Aut } _{{\cal D}_Q } (Q)
\right\} $.

\end{lem}

{\bf Proof.}
First note that,
for every $\Delta \in
{\rm Aut } _{{\cal D}_Q } (Q)$, the fact that
$\Delta $ fixes all sets $A_i ^j $ implies that
$\Delta ^U \in {\rm Aut } _{{\cal D}} (U)$ and that
$\Delta ^U _n ={\rm id} _{U_n } $.

Conversely, if
$\Psi \in {\rm Aut} _{\cal D} (U)$
satisfies $\Psi _n ={\rm id} _{U_n } $, then
$\Psi $ maps every $A_i ^j $ to itself.
Hence $\Psi |_Q \in {\rm Aut } _{{\cal D}_Q } (Q)$
and $\Psi =(\Psi |_Q )^U $.
\qed

\vspace{.1in}

So far, we have investigated the
structured ordered sets $(U,{\cal D})$,
$(U_n , {\cal D} _n )$, and
$(Q, {\cal D} _Q )$, and
groups directly associated with them.
Unfortunately,
${\rm Aut} _{{\cal D} _n } (U_n )$
can be a proper subset
of $\{ \Phi _n :  \Phi \in {\rm Aut} _{\cal D} (U)\} $.
Therefore, to
fully utilize the power of the ideas
presented so far for a recursive
deconstruction of ${\rm Aut} _{\cal D} (U)$, we
must focus on subgroups
$G^* $ of ${\rm Aut} _{\cal D} (U)$.

\begin{define}
\label{starautsetdef}

Let $G^* $ be a subgroup of ${\rm Aut} _{\cal D} (U)$.
We define
${\rm Aut } _{{\cal D}_n} ^* (U_n )
:=
\left\{ \Phi _n :\Phi \in G^* \right\} $,
${\rm Aut} _{{\cal D} _Q } ^{U*} (Q):=
\left\{ \Psi \in G^*: \Psi _n ={\rm id} _{U_n } \right\}
$
and
${\rm Aut} _{{\cal D} _Q } ^* (Q):=
\left\{ \Psi |_Q :\Psi \in {\rm Aut} _{{\cal D} _Q } ^{U*} (Q)\right\}
$.

\end{define}

\begin{theorem}
\label{pruneorbit6}

Let $G^* $ be a subgroup of ${\rm Aut} _{\cal D} (U)$.
The set ${\rm Aut} _{{\cal D} _Q } ^{U*} (Q)$ is a normal subgroup of
$G^* $
and
the factor group
$
G^* /{\rm Aut} _{{\cal D} _Q } ^{U*} (Q)$
is isomorphic to
${\rm Aut } _{{\cal D}_n} ^* (U_n )$.
Consequently
$
|G^*|
=
\left|{\rm Aut } _{{\cal D}_n} ^* (U_n )\right|
\left| {\rm Aut } _{{\cal D}_Q } ^* (Q)\right| $.

\end{theorem}

{\bf Proof.}
First note that, because
$\Phi _n \Psi _n =(\Phi \Psi )_n $ and
$(\Phi ^{-1} )_n =(\Phi _n )^{-1} $, we have that
${\rm Aut } _{{\cal D}_n} ^* (U_n )$ is a subgroup
of
$G^* $
and that
${\rm Aut} _{{\cal D} _Q } ^{U*} (Q)$ is a subgroup of
$G^* $.

Let
$\Delta ^U \in {\rm Aut} _{{\cal D} _Q } ^{U*} (Q)$
and let $\Phi \in
G^* $.
Then
$
\left( \Phi ^{-1} \Delta ^U \Phi \right) _n
=
\Phi ^{-1} _n \Delta ^U _n \Phi _n
=
\Phi ^{-1} _n \Phi _n
=
{\rm id} _{U_n } $.
Hence
$\Phi ^{-1} {\rm Aut} _{{\cal D} _Q } ^{U*} (Q) \Phi
=
{\rm Aut} _{{\cal D} _Q } ^{U*} (Q)$
and therefore ${\rm Aut} _{{\cal D} _Q } ^{U*} (Q)$
is a normal subgroup of
$G^* $.

Moreover, for all
$\Phi, \Psi \in {\rm Aut} _{\cal D} ^* (U)$, we have
$\Phi {\rm Aut} _{{\cal D} _Q } ^{U*} (Q)=
\Psi {\rm Aut} _{{\cal D} _Q } ^{U*} (Q)$
iff $\Phi _n =\Psi _n $.
Hence the factor group
$
G^*/{\rm Aut} _{{\cal D} _Q } ^{U*} (Q)$
is, via $\Phi {\rm Aut} _{{\cal D} _Q } ^{U*} (Q)\mapsto \Phi _n $,
isomorphic to
${\rm Aut } _{{\cal D}_n} ^* (U_n )$.
The
equation
now follows because
$\Psi \mapsto \Psi |_Q $ is an isomorphism from
${\rm Aut} _{{\cal D} _Q } ^{U*} (Q)$ to ${\rm Aut} _{{\cal D} _Q } ^{*} (Q)$.
\qed

\vspace{.1in}

The notations carry over
to subgroups $G^* $.

\begin{define}
\label{starSG}

Let $G^* $ be a subgroup of ${\rm Aut} _{\cal D} (U)$.
We define ${\cal D} _n ^* $ to be the set of
${\rm Aut } _{{\cal D}_n} ^* (U_n )$-orbits,
and
we define ${\cal D} _Q ^* $ to be the set of
${\rm Aut } _{{\cal D}_Q } ^* (Q)$-orbits.

\end{define}

\begin{define}

Let
$(P,{\cal D} )$ be a structured ordered set and let
$G^* $ be a subgroup of ${\rm Aut} _{\cal D} (P)$.
For all $S\subseteq P$, we define
$\Lambda _{\cal D} ^* (S):=\{ \Phi |_S :\Phi \in G^* \} $.

\end{define}

Because all we did is ``re-tighten" the interdependent orbit
unions (if this is even necessary), we have
${\rm Aut } _{{\cal D}_n ^* } ^* (U_n )=
{\rm Aut } _{{\cal D}_n} ^* (U_n )$
and
${\rm Aut } _{{\cal D}_Q ^* } ^* (Q)
={\rm Aut } _{{\cal D}_Q } ^* (Q)$.

\begin{define}

Let $G^* $ be a subgroup of ${\rm Aut} _{\cal D} (U)$.
The {\bf separation partition
${\cal S}^* (D_n )$
of $D_{n} $}
is the partition of $D_n $ that is contained in ${\cal D} _Q ^* $.

\end{define}

\begin{lem}
\label{pruneorbit5}

Let $G^* $ be a subgroup of ${\rm Aut} _{\cal D} (U)$.
Every $\Phi \in G^*$ respects
${\cal D} _Q ^* $.
Every nontrivial
${\cal D}_Q ^* $-orbit
$D$ that is contained
in
a set $A_i ^j $
is directly interdependent
with a
${\cal D}_Q ^* $-orbit
$S\in {\cal S} ^* (D_n )$, and all ${\cal D}_Q ^* $-orbits
that are directly interdependent with $D$ are contained in $D_n $.

\end{lem}

{\bf Proof.}
Suppose, for a contradiction, there are a
$\Phi \in G^* $
and an $S\in {\cal D} _Q ^* $ such that
$\Phi [S]\not\in {\cal D} _Q ^* $.
Then $\Phi [S]$ intersects two distinct sets $B,C\in {\cal D}_Q ^* $
or $\Phi [S]$ is strictly contained in a set $D\in {\cal D} _Q ^* $.
Because we are free to work with the inverse, it suffices to consider
the case
in which
$\Phi [S]$ intersects two distinct sets $B,C\in {\cal D} _Q ^* $.
Because $S\in {\cal D} _Q ^* $ is an ${\rm Aut} _{{\cal D} _Q ^*} ^* (Q)$-orbit,
there is a
$\Delta ^U \in {\rm Aut} _{{\cal D}_Q ^* } ^* (Q)$
that maps a $b\in \Phi ^{-1} [B]\cap S$ to a $c\in \Phi ^{-1} [C]\cap S$.
Now $\Phi (b)\in B$,
$\Phi \Delta ^U \Phi ^{-1} \left( \Phi (b)\right) =\Phi (c)\in C$
and $\Phi \Delta ^U \Phi ^{-1} |_Q \in {\rm Aut} _{{\cal D}_Q ^* } ^* (Q)$,
a contradiction to $B,C\in {\cal D}_Q ^* $.
We thus conclude that
every $\Phi \in G^* $ respects
${\cal D} _Q ^* $.

Finally,
by
Lemma \ref{pruneorbit3},
any nontrivial
${\cal D}_Q ^* $-orbit
that is contained
in
an $A_i ^j $
is directly interdependent
with an $S\in {\cal S}^* (D_n )$.
Because no two distinct $A_i ^j $ are directly interdependent and
no $A_i ^j $ is directly interdependent with a $D_k $ with $k>n$,
all orbits that are directly interdependent with an orbit in $A_i ^j $
must be contained in
$D_n $.
\qed

\section{Deconstructing Automorphisms}
\label{AuttoLambda}

Lemma \ref{expolem} below
is the motivation for the recursive deconstruction of
${\rm Aut} _{\cal D} (U)$ given in
Definition \ref{deconstruct}.
At every step, we will remove a noncutvertex of the orbit graph.
At the end of the deconstruction, we will be left with an interdependent orbit
union with $2$ orbits.
Because the removed vertices of the orbit graph are
not cutvertices,
the residual
structured ordered sets $(Q,{\cal D} _Q )$ have bipartite
orbit graphs.
Hence we start with
bipartite
orbit graphs.

\begin{lem}
\label{bipartiteOG}

Let $(P,{\cal D} )$ be a tight structured ordered set
such that every $D\in {\cal D}$
is directly interdependent with another
$E\in {\cal D}$, and
such that $P$ can be partitioned into sets $B$ and $T$ such that,
if $D,E\in {\cal D}$ are directly interdependent, then
$D\subseteq B$ and $E\subseteq T$
or
$D\subseteq T$ and $E\subseteq D$.
Let $G^* $ be a subgroup of
${\rm Aut} _{\cal D} (P)$.
Then
$\left| G^*\right|
=
\left| \Lambda _{\cal D} ^* (B)\right|
=
\left| \Lambda _{\cal D} ^* (T)\right|
$.

\end{lem}

{\bf Proof.}
For every $p\in P$, let $D_p $ be the unique
orbit in ${\cal D}$ that contains $p$.
If $p\in B$, then, because $p\in D_p \cap B$,
by hypothesis, we have $D_p \subseteq B$.
Similarly, if $p\in T$, then $D_p \subseteq T$.
We define a new order $\sqsubseteq $ on $P$ by,
for $b\in B$ and $t\in T$, setting
$b\sqsubseteq t$ iff
$D_b \updownharpoons _{\cal D} D_t $ and
$b\sim t$.

Let
$b_1 , b_2 \in B$ be two distinct elements
such that $D_{b_1} =D_{b_2} $.
To prove that $\{ b_1 , b_2 \} $ is not a
$\sqsubseteq $-order-autonomous
antichain, we argue as follows.
Because $(P,{\cal D})$ is tight,
and hence without slack,
without loss of generality,
there is an $a\in P\setminus \{ b_1 , b_2 \} $
such that $a\sim b_1 $ and $a\not\sim b_2 $.
In particular, this means that
$D_a \not= D_{b_1} =D_{b_2} $
satisfies
$D_a \updownharpoons _{\cal D} D_{b_1} =D_{b_2} $,
and hence
$a\sqsupseteq b_1 $
and
$a\not\sqsupseteq b_2 $.
Thus $\{ b_1 , b_2 \} $ is not a
$\sqsubseteq $-order-autonomous
antichain.
Similarly, for any
two distinct elements $t_1 , t_2 \in T$
such that $D_{t_1} =D_{t_2} $,
we have that $\{ t_1 , t_2 \} $ is not a
$\sqsubseteq $-order-autonomous
antichain.

Now let $\Phi \in {\rm Aut} _{\cal D} (P)$
and let $b\sqsubset t$.
Then
$b\in D_b \updownharpoons _{\cal D} D_t \ni t$ and $b\sim t$.
Because
$\Phi \in {\rm Aut} _{\cal D} (P)$, we have
$\Phi (b)\in D_b \updownharpoons _{\cal D} D_t \ni \Phi (t)$ and
$\Phi (b)\sim \Phi (t)$.
Hence
$\Phi (b)\sqsubset \Phi (t)$.
Because $\Phi $ is bijective,
$\Phi $ is a $\sqsubseteq $-automorphism.
In particular, every function in $G^* $
is a $\sqsubseteq $-automorphism.

Finally let $\Phi \in G^* $.
For any $t\in T$, we
apply
Lemma \ref{allbutoneac}
to $B\cup D_t $
and $\Phi _{B\cup D_t } $
with $A=D_t $ to conclude that
$\Phi |_{D_t } $ is determined
by $\Phi |_B $.
Thus every $\Phi \in G^* $ is uniquely determined
by $\Phi |_B $.
Similarly,
every $\Phi \in G^* $ is uniquely determined
by $\Phi |_T $.
We conclude that
$\left| G^*\right|
=
\left| \Lambda _{\cal D} ^* (B)\right|
=
\left| \Lambda _{\cal D} ^* (T)\right|
$.
\qed

\begin{lem}
\label{expolem}

Let $(U,{\cal D} )$ be a tight interdependent orbit union,
let $G^* $ be a subgroup of
${\rm Aut} _{{\cal D} } (U )$,
and let $D_m $ not be a cutvertex of ${\cal O} (U, {\cal D})$.
Let $b,c>0$ be so that the following hold.
\begin{enumerate}
\item
\label{expolem1}
$|{\rm Aut} _{{\cal D} _m } ^* (U_m )|\leq 2^{{b\over 2b-1}c|U_m |} $.
\item
\label{expolem2}
For all
$j\in \{ t,
\ldots , m\} $, we have
$|\Lambda _{{\cal D}_Q} ^* (D_j )|\leq 2^{c|D_j |} $.

\item
\label{expolem3}
For all $j\in \{ t, \ldots , m-1\} $
and $i\in \{ 1, \ldots , \ell _j \} $,
we have $\left| A_i ^j \right| \geq b$.
\end{enumerate}
Then
$\left| G^* \right| \leq 2^{{b\over 2b-1}c|U |} $.
If, in addition, we have $c\geq 1$
and $|{\rm Aut} _{{\cal D} _m } ^* (U_m )|\leq 2^{{4\over 7}c|U_m |} $,
then
$\left| G^* \right| \leq 2^{\min \left\{ {4\over 7} , {b\over 2b-1}\right\} c|U |} $.

\end{lem}

{\bf Proof.}
In case
$\left| {\rm Aut } _{{\cal D}_Q } ^* (Q)\right| =1$,
the
conclusions trivially follow from
Theorem \ref{pruneorbit6}.
For the
remainder, we can assume
that
$\left| {\rm Aut } _{{\cal D}_Q } ^* (Q)\right| >1$.
In particular, this means that $t<m$.

By Lemma \ref{pruneorbit5},
no two distinct ${\cal D}_Q $-orbits
in $\bigcup _{j=t} ^{m-1} D_j $ are directly interdependent
and every ${\cal D}_Q $-orbit
in $\bigcup _{j=t} ^{m-1} D_j $ is directly interdependent with
a ${\cal D}_Q $-orbit in $D_m $.
Let
$D\subseteq D_m $ be a ${\cal D} _Q $-orbit.
Clearly, $D$ is not directly interdependent with any
${\cal D}_Q $-orbits contained in $D_m $.
Hence there is an $a\in \bigcup _{j=t} ^{m-1} D_j $
that is comparable to some,
but not all, elements
of $D$,
and $D$ is directly interdependent with a
${\cal D} _Q $-orbit in $\bigcup _{j=t} ^{m-1} D_j $.
Therefore, we can apply
Lemma \ref{bipartiteOG} to $(Q, {\cal D} _Q )$
with $B:=D_m $ and $T:=\bigcup _{j=t} ^{m-1} D_j $.
We obtain
$
\left| {\rm Aut } _{{\cal D}_Q } ^* (Q)\right|
=
\left| \Lambda _{{\cal D}_Q } ^* (B)\right|
=
\left| \Lambda _{{\cal D}_Q } ^* (T)\right|
=
\left| \Lambda _{{\cal D}_Q } ^* (D_m )\right|
=
\left| \Lambda _{{\cal D}_Q } ^*
\left( \bigcup _{j=t} ^{m-1} D_j\right) \right|
$.

We first prove that
$\left| G^* \right| \leq 2^{{b\over 2b-1}c|U |} $.

By hypothesis
\ref{expolem2},
with $n_1 :=|D_m |$, we have
$\left| \Lambda _{{\cal D}_Q } ^* (D_m )\right| \leq 2^{cn_1 } $,
and with $n_2 :=\sum _{j=t} ^{m-1} |D_j |$, we
have
$\left| \Lambda _{{\cal D}_Q } ^*
\left( \bigcup _{j=t} ^{m-1} D_j\right) \right|
\leq
\prod _{j=t} ^{m-1}
\left| \Lambda _{{\cal D}_Q } ^*
\left( D_j\right) \right|
\leq
\prod _{j=t} ^{m-1}
2^{c|D_j |}
\leq
2^{cn_2 } $.
Thus
$
\left| {\rm Aut } _{{\cal D}_Q } ^* (Q)\right|
\leq
2^{c\min \{ n_1 , n_2 \} } $.

By hypothesis
\ref{expolem3},
for $j=t, \ldots , m-1$
and $i=1, \ldots , \ell _j $,
we have $\left| A_i ^j \right| \geq b$.
Hence $|U_m |\leq |U|-n_1 -n_2 +{1\over b} \max \{ n_1 , n_2 \} $.
Now we obtain the following
via Theorem \ref{pruneorbit6} and
hypothesis
\ref{expolem1}.
Because all computations are in the exponents,
we consider the base $2$ logarithm.
\begin{eqnarray*}
\lg \left( \left| G^* \right| \right)
& = &
\lg \left( \left| {\rm Aut } _{{\cal D}_m} ^* (U_m ) \right|
\left| {\rm Aut } _{{\cal D}_Q } ^* (Q)\right| \right)
\\
& = &
\lg \left( \left| {\rm Aut } _{{\cal D}_m} ^* (U_m ) \right| \right)
+
\lg \left( \left| {\rm Aut } _{{\cal D}_Q } ^* (Q)\right| \right)
\\
& \leq &
{b\over 2b-1} c|U_m |+c\min \{ n_1 , n_2 \}
\\
& \leq &
{b\over 2b-1} c\left( |U|-n_1 -n_2 +{1\over b} \max \{ n_1 , n_2 \} \right)
+
c(n_1 +n_2 -\max \{ n_1 , n_2 \} )
\\
& = &
{b\over 2b-1} c|U| +c\left( {b-1\over 2b-1} (n_1 +n_2 ) -
{2b-2\over 2b-1} \max \{ n_1 , n_2 \} \right)
\\
& \leq &
{b\over 2b-1} c|U|
\end{eqnarray*}

It remains to be proved that,
in case $c\geq 1$
and $|{\rm Aut} _{{\cal D} _m } ^* (U_m )|\leq 2^{{4\over 7}c|U_m |} $,
we have
$\left| G^* \right| \leq 2^{\min \left\{ {4\over 7} , {b\over 2b-1}\right\} c|U |} $.
Clearly, in this case, we can assume that $b\in \{ 2,3\} $.
Let
$Z$ be the number of
elements in $2$-antichains $A_i ^j $,
let $T$ be the number of
elements in $3$-antichains $A_i ^j $,
let $n_1 :=|D_m |$,
and let $n_2 :=\sum _{j=t} ^{m-1} |D_j |$.
Because every point in
$\bigcup _{j=t} ^{m-1} D_j\setminus (Z\cup T)$
is in a set $A_i ^j $ with at least $4$ elements, we obtain
\begin{eqnarray*}
|U_m |
& \leq &
|U|-n_1 -n_2 +{1\over 2} Z+{1\over 3} T +{1\over 4} \left( n_2 -Z-T\right)
\\
& = &
|U|-n_1 -{3\over 4}n_2
+{1\over 4} Z +{1\over 12} T.
\end{eqnarray*}

{\em Case 1:
$\displaystyle{n_2
\leq
n_1
+
{6c-{7\over 2}\over 4c} Z
+{20c- {56\over 3}\over 12c} T
}$.
}
Let
${\cal D} _Q ^2 $
and
${\cal D} _Q ^3 $
be comprised of
the sets $D_j $
that are partitioned
into sets $A_i ^j $ that contain $2$ or $3$ elements, respectively.
Let
${\cal D} _Q ^4 $
be comprised of
the sets
$D_j $
that are partitioned
into sets $A_i ^j $ that contain $4$ or more elements.
We obtain the following via Lemmas \ref{pruneorbit5} and \ref{bipartiteOG}.
\begin{eqnarray*}
\left| {\rm Aut } _{{\cal D}_Q } ^* (Q)\right|
& \leq &
\left| \Lambda _{{\cal D}_Q} ^* \left( \bigcup _{j=t} ^{m-1} D_j \right) \right|
\\
& \leq &
\prod _{D\in {\cal D} _Q ^2 } \left| \Lambda _{{\cal D}_Q} ^* (D )\right|
\prod _{D\in {\cal D} _Q ^3 } \left| \Lambda _{{\cal D}_Q} ^* (D )\right|
\prod _{D\in {\cal D} _Q ^4 } \left| \Lambda _{{\cal D}_Q} ^* (D )\right|
\\
& \leq &
2^{{1\over 2} Z} 6^{{1\over 3} T} 2^{c(n_2 -Z-T)}
\\
& \leq &
2^{{1\over 2} Z} 2^{{8\over 9} T} 2^{c(n_2 -Z-T)}
\\
& = &
2^{cn_2 +\left( {1\over 2}-c\right) Z
+\left( {8\over 9}-c\right) T}
.
\end{eqnarray*}

The above leads to the following estimate for
$|G^* |$.
Because all computations are in the exponent,
we work with the base 2 logarithm.
\begin{eqnarray*}
\lg \left( \left| G^* \right| \right)
& = &
\lg \left( \left| {\rm Aut } _{{\cal D}_m} ^* (U_m ) \right|
\left| {\rm Aut } _{{\cal D}_Q } ^* (Q)\right| \right)
\\
& = &
\lg \left( \left| {\rm Aut } _{{\cal D}_m} ^* (U_m ) \right| \right)
+
\lg \left( \left| {\rm Aut } _{{\cal D}_Q } ^* (Q)\right| \right)
\\
& \leq &
{4\over 7} c|U_m |
+
\left(
cn_2 +\left( {1\over 2}-c\right) Z
+\left( {8\over 9}-c\right) T \right)
\\
& \leq &
{4\over 7} c\left( |U|-n_1 -{3\over 4}n_2
+{1\over 4} Z +{1\over 12} T
\right)
+
{4\over 7} c
\left(
{7\over 4} n_2
+{7\over 4}
{{1\over 2}-c\over c} Z
+{7\over 4} { {8\over 9}-c\over c} T
\right)
\\
& = &
{4\over 7} c\left( |U|
-n_1 + n_2
+{1\over 4} Z
+{7\over 4}
{{1\over 2}-c\over c} Z
+{1\over 12} T
+{7\over 4} { {8\over 9}-c\over c} T
\right)
\\
& = &
{4\over 7} c\left( |U|
-n_1 + n_2
+{1\over 4c} cZ
+
{{7\over 2}-7c\over 4c} Z
+{1\over 12c} cT
+{ {56\over 3}-21c\over 12c} T
\right)
\\
& = &
{4\over 7} c\left( |U|
-n_1 + n_2
+
{{7\over 2}-6c\over 4c} Z
+{ {56\over 3}-20c\over 12c} T
\right)
\\
& = &
{4\over 7} c\left( |U|
+ n_2
-
\left(n_1
-
{{7\over 2}-6c\over 4c} Z
-{ {56\over 3}-20c\over 12c} T
\right)
\right)
\\
& = &
{4\over 7} c\left( |U|
+ n_2
-
\left(n_1
+
{6c-{7\over 2}\over 4c} Z
+{ 20c-{56\over 3}\over 12c}
\right)
\right)
\\
& \leq &
{4\over 7} c|U|,
\end{eqnarray*}
where the last step follows from the hypothesis for this case.

{\em Case 2:
$\displaystyle{n_2
>
n_1
+
{6c-{7\over 2}\over 4c} Z
+{20c- {56\over 3}\over 12c} T
}$.}
In this case,
by Lemma \ref{bipartiteOG},
we have $\left| {\rm Aut } _{{\cal D}_Q } ^* (Q)\right| \leq 2^{cn_1 } $,
by assumption
we have
$\displaystyle{
n_1-n_2 <
{{7\over 2}-6c\over 4c} Z
+{ {56\over 3}-20c\over 12c} T }$, and we obtain the following.
\begin{eqnarray*}
\lg \left( \left| G^* \right| \right)
& = &
\lg \left( \left| {\rm Aut } _{{\cal D}_m} ^* (U_m ) \right| \right)
+
\lg \left( \left| {\rm Aut } _{{\cal D}_Q } ^* (Q)\right| \right)
\\
& \leq &
{4\over 7} c|U_m |
+cn_1
\\
& \leq &
{4\over 7} c\left( |U|-n_1 -{3\over 4}n_2
+{1\over 4} Z +{1\over 12} T
\right)
+
{4\over 7}c\left( {7\over 4}n_1 \right)
\\
& = &
{4\over 7} c\left( |U|+{3\over 4}n_1 -{3\over 4}n_2
+{1\over 4} Z +{1\over 12} T
\right)
\\
& = &
{4\over 7} c\left( |U|+{3\over 4}(n_1 -n_2 )
+{1\over 4} Z +{1\over 12} T
\right)
\\
& \leq &
{4\over 7} c\left( |U|+{3\over 4}
\left(
{{7\over 2}-6c\over 4c} Z
+{ {56\over 3}-20c\over 12c} T
\right)
+{1\over 4} Z +{1\over 12} T
\right)
\\
& = &
{4\over 7} c\left( |U|+
{{21\over 2}-18c\over 16c} Z
+{ 56-60c\over 48c} T
+{4c\over 16c} Z +{4c\over 48c} T
\right)
\\
& = &
{4\over 7} c\left( |U|+
{{21\over 2}-14c\over 16c} Z
+{ 56-56c\over 48c} T
\right)
\qquad ({\rm recall \ } c\geq 1)
\\
& \leq &
{4\over 7} c|U|.
\end{eqnarray*}
\qed

\vspace{.1in}

We conclude this section with
some terminology to
facilitate the recursive application of Lemma
\ref{expolem}.

\begin{define}
\label{deconstruct}

Let $(U,{\cal D} )$ be a tight interdependent orbit union.
We define $(U^{m+1} , {\cal D} ^{m+1} ):=(U , {\cal D} )$.
Recursively,
for $j=m, \ldots , 3$, we can choose
$D^{j} $ to be a noncutvertex of
$\left( U^{j+1} , {\cal D} ^{j+1} \right) $.
We can define
$\left( U^{j} , {\cal D} ^{j} \right) $ to be the
pruned and compacted ordered set
from applying Definition \ref{phindef}
to $\left( U^{j+1} , {\cal D} ^{j+1} \right) $ and
$D^j $.
We can define
$\left( Q^{j} , {\cal D} _{Q^{j} } \right) $
to be the set
$(Q,{\cal D} _Q )$
obtained by applying
Definition \ref{setQdef}
to $\left( U^{j+1} , {\cal D} ^{j+1} \right) $ and
$D^j $.
A sequence
$\left\{ \left( U^j , {\cal D} ^j \right) \right\} _{j=m+1} ^3 $
obtained in this fashion is called
a {\bf deconstruction sequence} of $(U,{\cal D} )$.
The corresponding sequence
$\left\{ \left( Q^{j} , {\cal D} _{Q^{j} } \right) \right\} _{j=m} ^3 $
obtained in this fashion is called
the {\bf sequence of residuals}
for
$\left\{ \left( U^j , {\cal D} ^j \right) \right\} _{j=m} ^3 $.
Finally, we also call
$(U^3 , {\cal D} ^3 )$ the {\bf final residual}.

\end{define}

\begin{define}
\label{bdecon}

Let
$(U,{\cal D} )$
be a tight interdependent orbit union.
For $b>0$, a deconstruction sequence
$\left\{ \left( U^j , {\cal D} ^j \right) \right\} _{j=m+1} ^3 $
of $(U,{\cal D} )$
is called a {\bf $b$-deconstruction sequence} of $(U,{\cal D} )$
iff, for all $j=m+1 , \ldots , 3$,
all
nontrivial order-autonomous antichains
of $U^{j+1} \setminus D^j $
(also see Notation \ref{blocksinDj})
have at least $b$ elements.

\end{define}

Clearly, any tight interdependent orbit union has a
$2$-deconstruction sequence.

\begin{remark}
\label{useQinfuture}

{\rm
It should be noted that we have not used many details on the
the combinatorial structure of
the residuals $(Q^j ,{\cal D} _{Q^j} )$. An earlier draft of this paper, which
is now completely superseded by the application of
results on permutation groups in the next section,
extensively analyzed the structure of residuals
$(Q,{\cal D} _Q )$ to arrive at a, rather technical, proof of
Lemma \ref{orbits11orless} below.
Similar combinatorial work should allow for further advances on the
Automorphism Conjecture in the future.

It is easy to see that, when the removed orbit $D_m $ is
a noncutvertex in
${\cal O} (U, {\cal D} )$,
any two
interdependent orbit unions in $(Q,{\cal D} _Q )$ are
isomorphic.
Via Lemma \ref{bipartiteOG} the analysis
of
interdependent orbit unions in $(Q,{\cal D} _Q )$
can therefore be
reduced to
analyzing ordered sets of height $1$ such that, without loss of generality,
for any two orbits $D_1 $ and $D_2 $ that consist of minimal elements, the
groups $\Lambda _{\cal D} (D_1 )$
and $\Lambda _{\cal D} (D_2 )$ are permutation equivalent, that is, there
is an isomorphism that is induced by a bijection between the points of the
respective base sets $D_1 $ and $D_2 $.

}

\end{remark}

\section{Primitive Nestings}
\label{primnestsec}

We can now discuss some deep results from the theory of permutation groups
that should lead to significant advances on the Automorphism Problem.
In the case of Theorem \ref{marotiexpbound}, we present a
simple improvement that can be obtained with a rather trivial computation.
Theorem \ref{24special} rests on the shoulders of the giant which is the
Classification of Finite Simple Groups.
Consequently, the subsequent
results on the Automorphism Problem
are rather deep mathematics.
Thanks to the accessible writing in \cite{BuekLee}
and \cite{Maroti}, this depth is
within reach of the author's humble combinatorial means.

We start by reviewing a few key concepts
from the theory of permutation groups.

\begin{define}

Let $G$ be a permutation group on the set $X$.
We call $G$ {\bf transitive}
iff for all $x,u\in X$, there is a $\sigma \in G$ such that
$\sigma (x)=u$.

\end{define}

\begin{define}

Let $G$ be a permutation group on the set $X$.
We call $|X|$ the {\bf degree} of $G$.

\end{define}

Clearly, for a tight interdependent orbit union
$(U,{\cal D})$ and $D\in {\cal D}$, we have that
$\Lambda _{\cal D} (D)$ is a transitive permutation
group of degree $|D|$.

\begin{define}
\label{blockdef}

Let $G$ be a permutation group on the set $X$.
A subset $B\subseteq X$ is called a {\bf ($G$-)block}
iff, for all $\sigma \in G$, we have that
$\sigma [B]=B$ or $\sigma [B]\cap B=\emptyset $.
A block is called {\bf nontrivial} iff it is not
a singleton and not equal to
$X$.

\end{define}

By Lemma \ref{pruneorbit1}, for every
$j\in \{ s, \ldots , n-1\} $, the partition
${\cal A} ^j $ is a partition of $D_j $ into blocks of
$\Lambda _{\cal D} (D_j )$.
Let $G^* $ be a subgroup of ${\rm Aut} _{\cal D} (U)$.
By Lemma \ref{pruneorbit5}, for every
$j\in \{ t, \ldots , m\} $, the partition
of $D_j $ induced by ${\cal D} _Q ^* $
is a partition of $D_j $ into blocks of
$\Lambda _{{\cal D} _Q ^* } ^* (D_j )$.

\begin{define}
\label{primpg}

A permutation group $G$ is called {\bf primitive}
iff it has no nontrivial $G$-blocks.
A permutation group that is not
primitive is called
{\bf imprimitive}.

\end{define}

\begin{theorem}
\label{24special}

(See Corollary 1.4 in \cite{Maroti}.)
Let $G$ be a primitive permutation group of degree $n$ not containing
$A_n $. If $|G|>2^{n-1} $, then $G$ has degree at most $24$, and is permutation
isomorphic to one of the $24$ groups
in Table \ref{badgrouptable}
with their natural permutation
representation if not indicated otherwise in the table.

\end{theorem}

\begin{table}

\begin{tabular}{|l|r|r|r|r|r|}
\hline
Degree $n$	&	Group $G$	&	Order $|G|$	&	$\lg (|G|)/n\leq $	&	transitivity	\\ 	\hline \hline
5	&	AGL(1,5)	&	20	&	0.8644	&	s2	\\	\hline
6	&	PSL(2,5)	&	60	&	0.9845	&	2p	\\	
6	&	PGL(2,5)	&	120	&	1.1512	&	s3	\\	\hline
7	&	PSL(3,2)	&	168	&	1.0561	&	2	\\	\hline
8	&	A$\Gamma $L(1,8)	&	168	&	0.9241	&	2p	\\	
8	&	PSL(2,7)	&	168	&	0.9241	&	2p	\\	
8	&	PGL(2,7)	&	336	&	1.0491	&	s3	\\	
8	&	ASL(3,2)=AGL(3,2)	&	1344	&	1.2991	&	3	\\	\hline
9	&	AGL(2,3)	&	432	&	0.9728	&	2	\\	
9	&	PSL(2,8)	&	504	&	0.9975	&	s3p	\\	
9	&	P$\Gamma $L(2,8)	&	1512	&	1.1736	&	3p	\\	\hline
10	&	PGL(2,9)	&	720	&	0.9492	&	s3	\\	
10	&	$M_{10}$	&	720	&	0.9492	&	s3	\\	
10	&	$S_6 $ primitive on 10 elt.	&	720	&	0.9492	&	2p	\\	
10	&	P$\Gamma $L(2,9)	&	1440	&	1.0492	&	3	\\	\hline
11	&	$M_{11}$	&	7920	&	1.1774	&	s4	\\	\hline
12	&	$M_{11} $ on 12 elements	&	7920	&	1.0793	&	3p	\\	
12	&	$M_{12}$	&	95040	&	1.3781	&	s5	\\	\hline
13	&	PSL(3,3)	&	5616	&	0.9582	&	2	\\	\hline
15	&	PSL(4,2)	&	20160	&	0.9533	&	2	\\	\hline
16	&	$2^4 :A_7 $	&	40320	&	0.9563	&	3	\\	
16	&	$2^4 :L_4 (2) =$AGL(4,2)	&	322560	&	1.1438	&	3	\\	\hline
23	&	$M_{23}$	&	10200960	&	1.0123	&	4	\\	\hline
24	&	$M_{24}$	&	244823040	&	1.1612	&	5	\\	
\hline
\end{tabular}

\caption{The permutation groups
given in \cite{Maroti} which have degree $n$ and order
$\geq 2^{n-1} $.
Data from \cite{BuekLee}.
In case of mismatched names, both names
are given.
The upper bound for $\lg (|G|)/n$
is obtained by rounding up to the fourth digit.
The number $k$ in the transitivity column gives the
$k$-transitivity of the
group, ``s" stands for ``sharply," ``p" stands for ``primitive,"
and both are more stringent than $k$-transitivity.
}
\label{badgrouptable}

\end{table}

\begin{theorem}
\label{marotiexpbound}

(Compare with Corollary 1.2 in  \cite{Maroti}.)
If $G$ is a
primitive subgroup of $S_n $ not containing $A_n $,
then
$|G| < 2^{1.38n}$.
Moreover, if $n> 24$, then $|G| < 2^n $.

\end{theorem}

{\bf Proof.}
The only difference to
\cite{Maroti}, Corollary 1.2
is that
\cite{Maroti}, Corollary 1.2
guarantees, for
a
primitive subgroup $G$ of $S_n $ not containing $A_n $,
that
$|G| < 3^n $.
By Theorem \ref{24special}, we have that
$|G|\leq 2^{n-1} $
unless
$G$ is permutation
isomorphic to one of the $24$ groups
in Table \ref{badgrouptable}
with their natural permutation
representation if not indicated otherwise in the table.

For every group $G$
in Table \ref{badgrouptable},
the degree $n$ and
an upper bound $c_G $ for
${\lg (|G|)\over n} $
are given.
Clearly, $|G|\leq 2^{c_G n} $, and all values $c_G $ are
smaller than the factor $1.38$ used in the statement above.
\qed

\vspace{.1in}

For the following, the author apologizes for
possibly (even likely) restating some
items from the theory of permutation groups.
However,
standard works on permutation groups, such as \cite{Dixmort},
or references, such as \cite{Maroti},
formulate results
primarily in the language of groups,
which does not seem to have a simple connection to
interdependent orbit unions
$(U_m , {\cal D} _m ) $,
structured ordered sets $(Q,{\cal D} _Q )$,
or Lemma \ref{expolem}.
The following will
keep the presentation self contained, and,
in Lemmas \ref{1reductprimnestseq}
and \ref{1reductprimnestresid},
make the requisite connections.

\begin{define}

Let $G$ be a permutation group on the set $X$
and let $C\subseteq X$.
We define
$G\cdot C:=\{ \sigma [C]:\sigma \in G\} $ and
$G|_C :=\{ \sigma |_C:\sigma \in G, \sigma [C]\subseteq C\} $.

\end{define}

\begin{prop}
\label{folkloreonblocks}

(Folklore.)
Let $G$ be a
transitive
permutation group on the set $X$
and let $B$ be a $G$-block.
Then $G\cdot B$ is a partition of $X$,
every $\sigma \in G$ respects $G\cdot B$, and, for any two
$B_1 , B_2 \in G\cdot B$, we have that $G|_{B_1 } $
is isomorphic to $G|_{B_2 } $.
Finally, if $A$ is a block that contains $B$, then
every block in $G \cdot B$ that intersects $A$
is contained in $A$.

\end{prop}

{\bf Proof.}
Let $\sigma , \tau \in G$ such that
$\sigma [B]\cap \tau [B]\not= \emptyset $.
Then $\tau ^{-1} \sigma [B]\cap B\not= \emptyset $.
Hence, because $B$ is a block, we infer
$\tau ^{-1} \sigma [B]\cap B=B$, and then $\sigma [B]= \tau [B]$.
This proves that
$G\cdot B$ is a partition of $X$ and that
every $\sigma \in G$ respects $G\cdot B$.

Now let
$B_1 , B_2 \in G\cdot B$
and let
$\mu \in G$ be so that $\mu [B_1 ]=B_2 $.
The function $\Phi (\sigma ):=\mu \sigma \mu ^{-1} |_{B_2 } $
is a homomorphism from
the subgroup
$\{ \sigma \in G: \sigma [B_1 ]=B_1 \} $
onto
the group
$G|_{B_2} $.
Its kernel is the set
${\rm ker} (\Phi )=\{ \sigma \in G: \sigma |_{B_1 }= {\rm id} _{B_1 } \} $
and the Isomorphism Theorem now provides the claimed isomorphism.

Finally, let $A$ be a block that contains $B$
and let $C:=\sigma [B]$ intersect $A$.
Then $\sigma [A]\cap A\not= \emptyset $,
and therefore $\sigma [A]=A$.
Consequently, $C=\sigma [B]\subseteq \sigma [A]=A$.
\qed

\begin{define}

Let $G$ be a permutation group on the set $X$
and let $B\subseteq X$ be a block.
For $\sigma \in G$, we define
$\sigma \hochdamit ^{G\cdot B} $ by,
for all $C\in G\cdot B$, setting
$\sigma \hochdamit ^{G\cdot B} (C):=\sigma [C]$.
Let $A$ be a block that contains $B$.
We define
$A[G\cdot B]:=\{ C\in G\cdot B: C\subseteq A\} $
and
$G\hochdamit ^{G\cdot B} _A :=
\{ \sigma \hochdamit ^{G\cdot B} |_{A[G\cdot B]} :
\sigma \in G, \sigma [A]\subseteq A\} $.\footnote{Unless the
author is mitaken, the $G\hochdamit ^{G\cdot B} _A $ are sections of $G$.
However, the definition of sections might include other structures and,
to the author's knowledge, does not explicitly refer to the block structures
we define here to make the connection to Lemma \ref{expolem}.}

\end{define}

\begin{define}

Let $G$ be a
transitive permutation group
on the set $X$.
A sequence of
blocks
$\{ x\} =B_0 \subsetneq B_1 \subsetneq \cdots \subsetneq B_m = X$
such that,
for all $j\in \{ 0,\ldots , m-1\} $ we have that
$G\hochdamit ^{G\cdot B_j } _{B_{j+1} } $
is primitive
is called a {\bf primitive nesting}.

\end{define}

By starting with a singleton $B_0 =\{ x\} $
and successively choosing a smallest size block
that contains, but is not equal to, the block currently under consideration,
we see that every finite permutation group has indeed many primitive nestings.
Theorem \ref{marotiexpbound}
shows that, except for alternating or symmetric groups
$G\hochdamit ^{G\cdot B_j } _{B_{j+1} } $
and $G|_{B_1 } $, the sizes of these groups are bounded by
$2^{1.38n} $.
Similar to some results in
\cite{Maroti}, it makes sense to exclude
groups
$G\hochdamit ^{G\cdot B_j } _{B_{j+1} } $
that contain the alternating group
on their respective domains.

\begin{define}

Let $G$ be a
transitive permutation group
on the set $X$ and let
$\{ x\} =B_0 \subsetneq B_1 \subsetneq \cdots \subsetneq B_m = X$
be a primitive nesting for $G$.
If, for $j=0,\ldots ,m-1$, no $G\hochdamit ^{G\cdot B_j } _{B_{j+1} } $
contains a copy of
$A_{|B_{j+1} |\over |B_j |} (B_{j+1} [G\cdot B_j ])$ with
${|B_{j+1} |\over |B_j |}\geq 6$,
then the primitive nesting is called {\bf proper}.
Otherwise, the
primitive nesting is called {\bf improper}
and every index $j$ for which
$G\hochdamit ^{G\cdot B_j } _{B_{j+1} } $
contains a copy of $A_{|B_{j+1} |\over |B_j |} (B_{j+1} [G\cdot B_j ])$
is called an {\bf index of a factorial factor}.

\end{define}

Clearly, for a primitive permutation group
$G$ on the set $X$, the only primitive nestings are of the form
$\{ x\} \subsetneq X$, and they are proper iff $G$ does not contain
the alternating group $A_{|X|} $.
For imprimitive groups with a proper primitive nesting,
we record the following.

\begin{lem}
\label{controllednest}

Let $G$ be a
transitive permutation group
on the set $X$
such that
there
is a proper primitive nesting
$\{ x\} =B_0 \subseteq B_1 \subseteq \cdots \subseteq B_m = X$.
Then $|G|\leq 2^{1.7376 |X|} $.

\end{lem}

{\bf Proof.}
Because the nesting is proper, for every $j\in \{ 0,\ldots , m-1\} $,
we have
$\left| G\hochdamit ^{G\cdot B_j } _{B_{j+1} } \right|
\leq
2^{d_j {|B_{j+1} |\over |B_j |} } $, where $d_j $ is the maximum
factor for primitive groups of degree
${|B_{j+1} |\over |B_j |}
\in \{ 6,7,8,9,10,11,12,13,15,16,23,24\} $, which
can be found in Table \ref{badgrouptable} and is $\leq 1.38$,
or, for ${|B_{j+1} |\over |B_j |}\leq 5$,
$d_j :=\lg \left( \left( {|B_{j+1} |\over |B_j |}\right) !\right) /{|B_{j+1} |\over |B_j |}
\leq 1.3814$,
or, see Theorem \ref{24special}, for
degrees $\geq 14$ not covered so far,
$d_j :=1$.
Thus, for all $j$, we have
$d_j \leq 1.3814$.

By Proposition \ref{folkloreonblocks},
for any block $B$ and any two
$C_1 ,C_2 \in G\cdot B$, we have that $G|_{C_1 } $
is isomorphic to $G|_{C_2 } $.
Hence

\begin{eqnarray*}
|G|
& \leq &
\prod _{j=0} ^{m-1}
\left| G\hochdamit ^{G\cdot B_j } _{B_{j+1} } \right|
^{|X|\over |B_{j+1} |}
\\
& \leq &
\prod _{j=0} ^{m-1}
\left( 2^{d_j {|B_{j+1} |\over |B_j |} } \right)
^{|X|\over |B_{j+1} |}
\\
& = &
\prod _{j=0} ^{m-1}
2^{d_j {|X|\over |B_{j} |} }
\\
& = &
2^{|X|\left( \sum _{j=0} ^{m-1}
{d_j \over |B_{j} |} \right) }
\\
& = &
2^{|X|\left( \sum _{j=0} ^{s}
{d_j \over |B_{j} |} \right) }
2^{|X|\left( \sum _{j=s+1} ^{m-1}
{d_j \over |B_{j} |} \right) }
\\
& \leq &
2^{|X|\left( \sum _{j=0} ^{s}
{d_j \over |B_{j} |} \right) }
2^{|X|\left( \sum _{j=s+1} ^{m-1}
{1.38 \over 2^j } \right) }
\\
& \leq &
2^{|X|\left( \sum _{j=0} ^{s}
{d_j \over |B_{j} |} \right) }
2^{|X|1.38 \left( \sum _{j=s+1} ^{\infty }
{1\over 2^j } \right) }
\\
& = &
2^{|X|\left( \sum _{j=0} ^{s}
d_j \prod _{i=0} ^{j-1} {|B_i |\over |B_{i+1} |}  \right) }
2^{|X|1.38 {1\over 2^{s} }  }
\end{eqnarray*}

It is simple to write the nested loops that compute
an upper bound for
$\sum _{j=0} ^{s}
d_j \prod _{i=0} ^{j-1} {|B_i |\over |B_{i+1} |} $
for $s=7$:
For any degree
${|B_{j+1} |\over |B_j |}$ that is $\leq 5$ or tabulated
in Table \ref{badgrouptable},
we use $d_j $ as indicated above and we use
${|B_{i} |\over |B_{i+1} |}$ as a factor in the product.
For any degree
${|B_{i+1} |\over |B_i |}>5$ that is not tabulated
in Table \ref{badgrouptable},
we use $d_j =1$ and we use
${1\over 14} $ as an upper bound
for the factor
${|B_{i} |\over |B_{i+1} |}$
in the product.

We obtain
$\sum _{j=0} ^{7}
d_j \prod _{i=0} ^{j-1} {|B_i |\over |B_{i+1} |}
\leq 1.7268$.
Together with ${1.38\over 2^7 } \leq 0.0108$,
we obtain the claimed estimate.
\qed

\vspace{.1in}

Let $(U,{\cal D})$ be an interdependent orbit union.
As was mentioned after Definition \ref{blockdef},
by Lemma \ref{pruneorbit1}, for every
$j\in \{ s, \ldots , n-1\} $, the partition
${\cal A} ^j $ is a partition of $D_j $ into blocks of
$\Lambda _{\cal D} (D_j )$.
In a deconstruction sequence,
the same orbit $D\in {\cal D}$ can be split multiple
times in this fashion, and there may be further
ways to split $D$ into blocks of
$\Lambda _{\cal D} (D_j )$ which may or may not arise
in deconstruction sequences, and which might also not be nested
with blocks that arise as sets $A_i ^j $.
(Also see Remark \ref{refinekeylemmss} below on this subject.)

To connect Lemma \ref{controllednest}
with
Lemma \ref{expolem},
we must assure that primitive nestings split into
primitive nestings and that
no new primitive nestings could arise when
a $\Lambda _{\cal D} (D)$ is split by
a partition into sets $A_i ^j $.

\begin{lem}
\label{subgroupinheritsnesting}

Let $G$ be a
transitive permutation group
on the set $X$, let
$\{ x\} =B_0 \subsetneq B_1 \subsetneq \cdots \subsetneq B_m = X$
be a primitive nesting for $G$ and let
$k\in \{ 1, \ldots , m-1\} $.
For $j=k, \ldots , m$, define
$\widehat{B_j } :=\{ C\in G\cdot B_k : C\subseteq B_j \} $.
Then
$\{ B_k \} =\widehat{B_k } \subsetneq \widehat{B_{k+1} }
\subsetneq \cdots \subsetneq \widehat{B_m }
=G\cdot B_k $
is a primitive nesting for
$G\hochdamit _X ^{G\cdot B_k } $
and all primitive nestings of
$G\hochdamit _X ^{G\cdot B_k } $
are of this form.

\end{lem}

{\bf Proof.}
The fact that
$\{ B_k \} =\widehat{B_k } \subsetneq \widehat{B_{k+1} }
\subsetneq \cdots \subsetneq \widehat{B_m } =G\cdot B_k $
is a primitive nesting for
$G\hochdamit _X ^{G\cdot B_k } $
follows straight from the definitions, as
the only adjustment is that
$G\hochdamit _X ^{G\cdot B_k } $
acts on the blocks in
$G\cdot B_k $, whereas $G$ acts on the
elements of $X$.

Conversely, let
$\{ D_k \} =C_k \subsetneq C_{k+1}
\subsetneq \cdots \subsetneq
C_n  =G\cdot B_k $
be a primitive nesting
for
$G\hochdamit _X ^{G\cdot B_k } $.
For $j\in \{ k+1\, \ldots , n\} $, let
$D_j :=\bigcup C_j $.
We can now
find a primitive nesting
$\{ x\} =
E_0 \subsetneq E_1 \subsetneq \cdots \subsetneq E_q =D_k
\subsetneq D_{k+1} \subsetneq \cdots \subsetneq D_n $
for $G$.
(Because different block structures could be chosen,
we need not have $q=k$, but, as before, the rest
follows straight from the
definitions.)
This is the requisite primitive nesting for $G$ that induces
$\{ D_k \} =C_k \subsetneq C_{k+1}
\subsetneq \cdots \subsetneq
C_n  =G\cdot B_k $
as indicated.
\qed

\begin{lem}
\label{1reductprimnestseq}

Let $(U,{\cal D} =\{ D_1 , \ldots , D_m \} )$ be
a tight
interdependent orbit union such that $D_m $ is not a cutvertex of the
orbit graph ${\cal O} (U,{\cal D} )$.
Let $G^* $ be a
subgroup of
${\rm Aut} _{\cal D} (U)$
such that, for all $D\in {\cal D} $,
$\Lambda _{\cal D} ^* (D )$
is
transitive and
every primitive nesting of
$\Lambda _{\cal D} ^* (D )$
is proper.
Then,
${\cal D} _m ^* = {\cal D} _m $ and,
for all $D\in {\cal D} _m ^* $,
$\Lambda _{{\cal D}_m ^* } ^* (D )$
is transitive and
every primitive nesting for
$\Lambda _{{\cal D}_m ^* } ^* (D )$
is proper.

\end{lem}

{\bf Proof.}
We can assume that
${\cal D} =\{ D_1 , \ldots , D_m \} $
is labeled as in Notation \ref{standardOGnotation}.
We will use the notation from Section \ref{orbgraphsec}.
For $j\in \{ 1, \ldots , m-1\} $,
transitivity of the
$\Lambda _{{\cal D} } ^* (D_j )$ implies that
${\cal D} _m ^* ={\cal D} _m $.
Hence,
for all $D_j\cap U_m \in {\cal D} _m ^* $,
we have that
$\Lambda _{{\cal D}_m ^* } ^* (D_j\cap U_m )$
is transitive.

For $j<t$, we have that
$D_j = D_j\cap U_m $ and
$\Lambda _{{\cal D}_m ^* } ^* (D_j\cap U_m )=\Lambda _{{\cal D}} ^* (D_j )$,
which establishes that, for $j<t$,
all primitive nestings
of $\Lambda _{{\cal D}_m ^* } ^* (D_j )$
are proper.

Let
$j\in \{ t, \ldots , m-1\} $ and let
$\{ b \} =B_0 \subsetneq B_1 \subsetneq \cdots \subsetneq B_z
=D_j \cap U_m $
be a primitive nesting for
$\Lambda _{{\cal D}_m ^* } ^* (D_j \cap U_m )$.
For each $B_k $ define $\widehat{ B_k } :=
\{ A_i ^j :a_i ^j \in B_k \} $.
Let $\widetilde{B_0 } $ be the unique element of
$\widehat{B_0 } $.
Note that
the only change
in going from
$\{ b \} =B_0 \subsetneq B_1 \subsetneq \cdots \subsetneq B_z
=D_j \cap U_m $
to
$\left\{ \widetilde{B_0 } \right\} =\widehat{B_0 }\subsetneq
\widehat{B_1 }\subsetneq
\cdots \subsetneq \widehat{B_z }
={\cal A}^j
=\Lambda _{{\cal D}} ^* (D _j )\cdot \widetilde{B_0 } $
is that
permutations in
$\Lambda _{{\cal D}_m ^* } ^* (D_j \cap U_m )$
are traded for permutations of blocks in
$\Lambda _{{\cal D}} ^* (D_j )\hochdamit _{D_j }
^{\Lambda _{{\cal D}} ^* (D _j )\cdot \widetilde{B_0 } } $.
Hence
$\left\{ \widetilde{B_0 } \right\}
=\widehat{B_0 }\subsetneq \widehat{B_1 }\subsetneq
\cdots \subsetneq \widehat{B_z }
={\cal A}^j
=\Lambda _{{\cal D}} ^* (D _j )\cdot \widetilde{B_0 } $
is a primitive nesting for
$\Lambda _{{\cal D}} ^* (D_j )\hochdamit _{D_j }
^{\Lambda _{{\cal D}} ^* (D _j )\cdot \widetilde{B_0 } } $.
By Lemma \ref{subgroupinheritsnesting},
this primitive nesting is isomorphic to the
upper part of a primitive nesting for
$\Lambda _{{\cal D}} ^* (D _j )$, and hence it is proper.
Consequently
$\{ b \} =B_0 \subsetneq B_1 \subsetneq \cdots \subsetneq B_z
=D_j \cap U_m $ is proper, too.
\qed

\begin{lem}
\label{blockinheritsnesting}

Let $G$ be a
transitive permutation group
on the set $X$, let
$B$ be a block of $G$ and let
$\{ x\} =B_0 \subsetneq B_1 \subsetneq \cdots \subsetneq B_z = B$
be a primitive nesting for $G|_B$.
Then there is a primitive nesting
for $G$ that contains the above primitive nesting.
Moreover, if the containing primitive nesting is proper, then so is
$\{ x\} =B_0 \subsetneq B_1 \subsetneq \cdots \subsetneq B_z = B$.

\end{lem}

{\bf Proof.}
Any block of $G|_B $ is a block of $G$, too.
Thus
$B_0 \subsetneq B_1 \subsetneq \cdots \subsetneq B_z$
are blocks of $G$.
Conversely, any block of $G$ that is contained in $B$ is a block of
$G|_B $.
Therefore, there is no $j\in \{ 1, \ldots , z\} $
such that there is a block $C$ of $G$
such that $B_{j-1} \subsetneq C\subsetneq B_j $.
Now
$B_0 \subsetneq B_1 \subsetneq \cdots \subsetneq B_z$
can be extended to a primitive nesting for $G$.
Clearly, if this
primitive nesting is proper, then so is
$\{ x\} =B_0 \subsetneq B_1 \subsetneq \cdots \subsetneq B_z = B$.
\qed

\begin{lem}
\label{1reductprimnestresid}

Let $(U,{\cal D} )$ be
a tight
interdependent orbit union
such that $D_m $ is not a cutvertex of the
orbit graph ${\cal O} (U,{\cal D} )$.
Let $G^* $ be a subgroup of
${\rm Aut} _{\cal D} (U)$
such that, for all $D\in {\cal D} $,
$\Lambda _{\cal D} ^* (D )$
is
transitive and
every primitive nesting of
$\Lambda _{\cal D} ^* (D )$
is proper.
Then, for all $D\in {\cal D} _Q ^* $,
$\Lambda _{{\cal D}_Q ^* } ^* (D )$
is transitive and
every primitive nesting for
$\Lambda _{{\cal D}_Q ^* } ^* (D )$
is proper.
In particular,
for all
$j\in \{ t, \ldots , m-1\} $, we have
$|\Lambda _{{\cal D}_Q} ^* (D_j )|\leq 2^{1.7376|D_j |} $.

\end{lem}

{\bf Proof.}
Let $D\in {\cal D} _Q ^* $.
By definition, $D$ is an orbit of ${\rm Aut} _{{\cal D} _Q ^* } ^* (Q)$,
and hence
$\Lambda _{{\cal D}_Q ^* } ^* (D )$
acts transitively on $D$.
By Lemma \ref{pruneorbit5},
every $\Phi \in G^*$ respects
${\cal D} _Q ^* $, which implies that
$D$ is
a block
of $\Lambda _{{\cal D} }^* (D_j )$.

By Lemma \ref{blockinheritsnesting}, any primitive nesting
for
$\Lambda _{{\cal D}_Q ^* } ^* (D )$
is contained in a, necessarily proper, primitive nesting for
$\Lambda _{{\cal D} } ^* (D )$, and hence
it is proper.
Consequently, by Lemma \ref{controllednest},
for every
$D\in {\cal D} _Q ^* $, we have
$|\Lambda _{{\cal D}_Q} ^* (D )|\leq 2^{1.7376|D |} $.
Now,
for all
$j\in \{ t, \ldots , m-1\} $, we have
$
|\Lambda _{{\cal D}_Q} ^* (D_j )|
\leq
\prod _{D\in {\cal D}_Q ^* , D\subseteq D_j }
|\Lambda _{{\cal D}_Q} ^* (D )|
\leq
\prod _{D\in {\cal D}_Q ^* , D\subseteq D_j }
2^{1.7376|D |}
=
2^{1.7376|D_j |} $.
\qed

\begin{theorem}
\label{permgrtoordaut}

Let $(U,{\cal D} )$ be
a tight
interdependent orbit union
such that, for all $D\in {\cal D} $,
every primitive nesting for
$\Lambda _{\cal D} (D )$ is proper.
If $(U,{\cal D} )$ has a $b$-deconstruction sequence,
then
$
|{\rm Aut} _{\cal D} (U) |
\leq 2^{\min \left\{ {4\over 7}, {b\over 2b-1}\right\} 1.7376|U |} $.
In case $|{\cal D} |=2$, we even have
$
|{\rm Aut} _{\cal D} (U) |
\leq 2^{{1\over 2}1.7376|U |} $.

\end{theorem}

{\bf Proof.}
We will prove by induction on $m=|{\cal D} |$ that,
for every tight
interdependent orbit union $(U,{\cal D} )$
with a $b$-deconstruction sequence, and for
every subgroup
$G^* $ of ${\rm Aut} _{\cal D} (U)$
such that, for all $D\in {\cal D} $,
$\Lambda _{\cal D} ^* (D )$
is
transitive and
every primitive nesting of
$\Lambda _{\cal D} ^* (D )$
is proper,
we have
$
|G^* |
\leq 2^{\min \left\{ {4\over 7}, {b\over 2b-1}\right\} 1.7376|U |} $, and that,
in case $m=2$, we even have
$
|G^* |
\leq 2^{{1\over 2}1.7376|U |} $.

{\em Base step $m=2$:}
By
Lemma \ref{controllednest}, for $i=1,2$, we obtain
$\left| \Lambda _{{\cal D} } ^* (D_i )\right| \leq 2^{1.7376|D_i |} $.
Via Lemma \ref{bipartiteOG}, with $B=D_1 $, $T=D_2 $ and
$c=1.7376$, we obtain
$\left| G^*\right|
=
\left| \Lambda _{\cal D} ^* (D_1 )\right|
=
\left| \Lambda _{\cal D} ^* (D_2 )\right|
\leq
2^{1.7376\min \{ |D_1 |, |D_2 |\} }
\leq
2^{1.7376{1\over 2} |U| }
$.

{\em Induction step $(m-1)\to m$:}
Let ${\cal D} =\{ D_1 , \ldots , D_m \} $, labeled such that
$D_m $ is the
orbit that is removed in the first step of a
$b$-deconstruction sequence.
Then clearly,
$(U_m ,{\cal D}_m )$
has a
$b$-deconstruction sequence.
By Lemma \ref{1reductprimnestseq},
${\cal D} _m ^* = {\cal D} _m $ and,
for all $D\in {\cal D} _m ^* $,
$\Lambda _{{\cal D}_m ^* } ^* (D )$
is transitive and
every primitive nesting for
$\Lambda _{{\cal D}_m ^* } ^* (D )$
is proper.
Thus, by induction hypothesis,
$
|{\rm Aut} _{{\cal D}_m } ^* (U_m ) |
\leq 2^{\min \left\{ {4\over 7}, {b\over 2b-1}\right\} 1.7376|U_m |} $.

By Lemma \ref{1reductprimnestresid},
the claim follows from Lemma \ref{expolem}
with $c=1.7376$.
\qed

\vspace{.1in}

It is not surprising that bounds for
the automorphism group of the ordered set will be a bit tighter
than the bounds for the
permutation groups induced on the
individual orbits:
The direct interdependences of an orbit $D$ with other orbits
bring more points into play, and that is all that is used in
Theorem \ref{permgrtoordaut}.

\begin{remark}
\label{smallprimitivefactors}

{\rm
Because every deconstruction sequence is a 2-deconstruction sequence,
Theorem \ref{permgrtoordaut} guarantees
$
|{\rm Aut} _{\cal D} (U) |
\leq
2^{{4\over 7}1.747|U |}
\leq
2^{0.999|U |}
$ when all primitive nestings of all $\Lambda _{\cal D} (D)$
are proper.

Example (3) in \cite{DRSW} shows that there is a family of ordered sets with
$|{\rm End} (P)|=\left( \sqrt{4 +\sqrt{11} } \right) ^{|P|}
\leq 2^{1.4356|P|} $.
A lower bound on the number of endomorphisms of
the form $2^{|P|} $
would be very helpful indeed, and Problem (2) in \cite{DRSW}
asks if there is an infinite family of ordered sets with
fewer than $2^{1.2716|P|} $ endomorphisms.
However, until such bounds are established,
we cannot exclude that
there may be families of ordered sets with
fewer endomorphisms.

Hence, although
the
``exponential bound $2^{0.999|P|} $ unless there are blocks with 6
or more elements which carry
their alternating group" from
Theorem \ref{permgrtoordaut}
is a significant (and, given the nature of the results needed to
prove it, deep)
step forward for the theory of automorphisms of ordered sets,
further analysis is needed to
resolve the Automorphism Conjecture.
The challenge for better bounds for the number of
automorphisms lies in factorial factors
and, possibly,
in certain small primitive factors that could be
excluded to improve the estimate in Lemma \ref{controllednest}:
If we, in addition to $A_n $ and $S_n $ for $n\geq 6$,
also exclude
$A_4 $, $S_4 $, $A_5 $, $S_5 $,
$PGL(2,5)$, $AGL(3,2)$, and $M_{12} $ from
all primitive nestings, the same argument we gave here
would lead to a bound of $2^{0.76|U|} $.
This is less than the new
lower bound for the number of endomorphisms that will be established in
Theorem \ref{moreendos} below.

}

\end{remark}

\begin{remark}
\label{refinekeylemmss}

{\rm
The hypothesis in
Theorem \ref{permgrtoordaut}
that, for all $D\in {\cal D} $,
every primitive nesting for
$\Lambda _{\cal D} (D )$ is proper is rather strong,
but, via Lemmas \ref{1reductprimnestseq}
and \ref{1reductprimnestresid},
it facilitates the
recursive application of
Lemma \ref{expolem}.
The problem we face in the use of
Lemma \ref{expolem}
is that a block in a $\Lambda _{\cal D} (D_j )$
may be split by
${\cal A} ^j $ in a way that
leaves too large a factor for
${\rm Aut} _{{\cal D} _m } (U_m )$.

For example, consider the set
$X:=\left\{ \left( i,\{ j,k\} \right) : i,j,k\in \{ 1, \ldots , n\} , j\not= k\right\} $
with the permutation group $G$ that consists of
all permutations $\widehat{\sigma } (i,j):=(\sigma (i),\{ \sigma (j), \sigma (k)\} )$,
where $\sigma \in S_n $. Clearly, all rows
$R_i :=\{ (i,\{ j,k\} ): j,k\in \{ 1, \ldots , n\} , j\not= k\} $
and columns
$C_{\{ j,k\} } :=\{ (i,\{ j,k\} ): i\in \{ 1, \ldots , n\} \} $
are blocks of $G$.
Moreover, $|G|=n!$, $|X|=n\pmatrix{n\cr 2\cr } ={1\over 2} \left( n^3 -n^2 \right) $
and, as $n$ grows, $2^{{1\over 2} \left( n^3 -n^2 \right) }\gg n!$.
However, if $G$ were to be split by a partition ${\cal A} ^j $
that partitions into rows, then
we would have $n!$ permutations on the $n$ elements in the
corresponding orbit of the set $(U_m , {\cal D} _m )$.
A split into the columns would still
give $n!$ permutations on the $n$ elements in the
corresponding orbit, but we would have ${1\over 2} \left( n^2 -n\right) $
elements to facilitate an adequate estimate.
In either case, the corresponding groups on the residuals
consist of the identity.

It should be possible to address some of these issues with
a more detailed refinement of
Lemma \ref{expolem}
and of
the definitions
for deconstruction sequences
given here, but
this should be addressed in future work.
}

\end{remark}

\section{
The Automorphism Conjecture for Ordered Sets of Width 11}
\label{ACsmwidthproof}

The results presented so far provide substantial new insights into the
structure of automorphism groups of ordered sets.
Almost as a ``proof of concept,"
we now turn to confirming the Automorphism Conjecture for
ordered sets of width up to 11.

We will call an interdependent union
$(U, {\cal D} )$ {\bf max-locked}
iff $U$ is a max-locked ordered set.

\begin{lem}
\label{orbits11orlessmaxlocked}

Let $(U,{\cal D} )$ be an interdependent orbit union
of width $\leq 11$ such that there is a $D\in {\cal D}$ with
$|D|\geq 6$ such that
$\Lambda _{\cal D} (D)$ contains $A_{|D|} (D)$.
Then
$(U,{\cal D} )$ is max-locked.

\end{lem}

{\bf Proof.}
Let $D\in {\cal D}$ with
$|D|\geq 6$ such that
$\Lambda _{\cal D} (D)$ contains $A_{|D|} (D)$
and let
$E\updownharpoons _{\cal D} D$ be directly interdependent with
$D$.
By Lemma \ref{firstauffaecherung},
because $|E|\leq 11<15=\pmatrix{ 6\cr 2\cr } $,
we conclude, with
$w=|D|$ that $D\cup E$ is isomorphic to
$wC_2 $ or $S_w $.
In particular, $|E|=|D|$ and
$\Lambda _{\cal D} (E)$ contains $A_{|E|} (E)$.

Because the above applies to
any $D\in {\cal D}$ with
$|D|\geq 6$ such that
$\Lambda _{\cal D} (D)$ contains $A_{|D|} (D)$
and because ${\cal O} (U,{\cal D} )$ is connected,
we conclude that,
for any two $D,E\in {\cal D}$ with
$D\updownharpoons _{\cal D} E$, we have that
$D\cup E$ is isomorphic to
$wC_2 $ or $S_w $.
Finally, because any two orbits of size at least $6$ must have
comparable elements, we conclude that
$(U,{\cal D} )$ is max-locked.
\qed

\begin{lem}
\label{orbits11orless}

Let $(U,{\cal D} )$ be a
tight interdependent orbit union
such that, for all $D\in {\cal D}$, we have $|D|\leq 11$, and all
orbits $|D|$ such that $\Lambda _{\cal D} (D)$  contains $A_{|D|} (D)$
satisfy $|D|\leq 5$.
Then
$|{\rm Aut} _{{\cal D} } (U )|\leq 2^{0.993|U |} $.

\end{lem}

{\bf Proof.}
Let $D\in {\cal D}$ such that $\Lambda _{\cal D} (D)$
contains a factorial factor on $n$ elements.
Because $|D|\leq 11$,
we have that $n\geq 6$ would imply that
$\Lambda _{\cal D} (D)$  contains $A_{|D|} (D)$, which was excluded.
Therefore $n\leq 5$.
Hence every primitive nesting for
$\Lambda _{\cal D} (D)$ is proper.

Because $D\in {\cal D}$ was arbitrary, we can
apply
Theorem \ref{permgrtoordaut} with $b=2$
and we obtain
$
|{\rm Aut} _{\cal D} (U) |
\leq 2^{{4\over 7}\cdot 1.7376|U |} \leq 2^{0.993|U |} $.
\qed

\vspace{.1in}

To prove the Automorphism Conjecture for
ordered sets of width up to $11$ in Theorem \ref{ACw<=12}
below,
we first establish a
lower bound on the number of endomorphisms in Lemma \ref{2tonforbddwidth}.
Then,
from
Proposition \ref{usingEndD}
through
Proposition \ref{manyinmaxlock}, we show
that a relative abundance of
max-locked interdependent orbit unions of
height $1$
guaranteees that the Automorphism Conjecture holds.
Lemma \ref{wid12AClem} combines the
work so far into an upper bound for the number of automorphisms
when there are no nontrivial order-autonomous antichains, and
Proposition \ref{onelexiter} shows that
a single execution of the lexicographic sum construction,
as occurs for example when nontrivial order-autonomous antichains are
inserted, does not affect the status of the Automorphism Conjecture.

\begin{lem}
\label{2tonforbddwidth}

Let $w\in {\mat N}$ and $\varepsilon >0$.
There is an $N\in {\mat N}$ such that
every
ordered set $P$ of width $\leq w$ with $n:=|P|\geq N$ elements
has at least
$2^{\left( 1
-\varepsilon \right) n} $
endomorphisms.

\end{lem}

{\bf Proof.}
Let $P$ be an ordered set of height $h$ with $n$ elements.
The proof of Theorem 1 in \cite{DRSW} (on page 20 of \cite{DRSW})
shows that $P$ has
at least $2^{{h\over h+1} n} $
endomorphisms that are surjective onto a chain
of length $h$.

Let $N\in {\mat N}$ be so that
${(N/w)-1\over [(N/w)-1]+1}>1-\varepsilon $ and let
$P$ be an ordered set of width $w$ with $n\geq N$ elements.
Then the height $h$ of $P$ satisfies
$h\geq {n\over w}-1 \geq {N\over w} -1$.
Hence
$P$ has
at least $2^{{h\over h+1} n} \geq 2^{(1-\varepsilon )n} $
endomorphisms.
\qed

\vspace{.1in}

Similar to ${\rm Aut} _{\cal D} (P)$,
we define
${\rm End} _{\cal D} (P)$.

\begin{define}
\label{EndDdef}

Let $(P, {\cal D} )$ be
a structured
ordered set.
We define ${\rm End} _{\cal D} (P)$ to be the set of
order-preserving maps $f:P\to P$ such that,
for all $D\in {\cal D}$, we have
$f[D]\subseteq D$.

\end{define}

\begin{prop}
\label{usingEndD}

Let $P$ be an ordered set and let
${\cal U}$ denote the set of
nontrivial natural
interdependent orbit unions of $P$. Then
$${
|{\rm Aut} (P)|
\over
|{\rm End} (P)|
}
=
{
|{\rm Aut}_{\cal N} (P)|
\over
|{\rm End}_{\cal N} (P)|
}
\leq
\prod _{U\in {\cal U} }
{
|{\rm Aut}_{{\cal N}|U} (U)|
\over
|{\rm End}_{{\cal N}|U} (U)|
}
.$$

\end{prop}

{\bf Proof.}
The functions in each
set ${\rm End}_{{\cal N}|U} (U)$ can be combined into
endomorphisms of $P$ in the same way as the functions in
${\rm Aut}_{{\cal N}|U} (U)$ are combined in
Proposition \ref{getallfromiou}.
\qed

\begin{lem}
\label{maxlockh1presrank}

$|{\rm End}_{\cal N} (wC_2 )|\geq w^w $
and
$|{\rm End}_{\cal N} (S_w )|\geq (w-1)^w $.

\end{lem}

{\bf Proof.}
The result about
$|{\rm End}_{\cal N} (wC_2 )|$ follows from the fact that every self map of
the minimal elements can be extended to a function in
${\rm End}_{\cal N} (wC_2 )$.
For the result about
$|{\rm End}_{\cal N} (S_w )|$, let $t\in S_w $ be maximal in $S_w $.
Now we can map the maximal elements of $S_w $
to $t$, we can map
every one of
the $w$ minimal elements of $S_w $
to any of the $w-1$ minimal elements in
$\downarrow t$, and, through each such choice,
we obtain a function in
${\rm End}_{\cal N} (S_w )$.
\qed

\begin{prop}
\label{manyinmaxlock}

The Automorphism Conjecture is true for the class
of ordered sets $P$
such that at least
$\lg (|P|)$ elements of $P$ are contained in
max-locked interdependent orbit unions of height $1$.

\end{prop}

{\bf Proof.}
Let $P$ be an ordered set
with $|P|=n$ elements
such that at least
$\lg (n)$ elements of $P$ are contained in
max-locked interdependent orbit unions of height $1$.

In case
$\geq \sqrt{\lg (n)} $ elements of $P$ are contained in
one max-locked interdependent orbit union
$U$ of height
$1$ and width $w$, we have
$w\geq {1\over 2}
\sqrt{\lg (n)} $.
By
Lemma \ref{maxlockh1presrank},
we obtain
$|{\rm End}_{{\cal N}|U} (U)|
\geq (w-1)^w $.
Moreover, $U$ has
exactly $w!$
automorphisms. By Proposition \ref{usingEndD}, we obtain the following.
$${ |{\rm Aut} (P)|
\over
|{\rm End} (P)|}
\leq
{w!\over (w-1)^w }
\leq
{\left( {1\over 2}
\sqrt{\lg (n)}\right) !\over
\left( {1\over 2}
\sqrt{\lg (n)}-1\right) ^{{1\over 2} \sqrt{\lg (n)}} }
,$$
and the right hand side goes to zero as $n\to \infty $.

In case no max-locked interdependent orbit union
$U\subseteq P$
of height
$1$
contains
$\geq \sqrt{\lg (n)} $ elements, we have that
$P$ contains
$\geq \sqrt{\lg (n)} $
max-locked
interdependent orbit unions
$(U, {\cal N}|U )$
of height $1$.
By Lemma \ref{maxlockh1presrank},
we have
${
|{\rm Aut}_{{\cal N}|U} (U)|
\over
|{\rm End}_{{\cal N}|U} (U)|
}
\leq
{w!\over (w-1)^w }
, $
which, for
$w\geq 3$,
is $\leq {3\over 4}$.
For
width $w=2$, we have
${
|{\rm Aut}_{{\cal N}|U} (U)|
\over
|{\rm End}_{{\cal N}|U} (U)|
}
\leq
{2\over 4}<{3\over 4} $.
By Proposition \ref{usingEndD},
we obtain
$${ |{\rm Aut} (P)|
\over
|{\rm End} (P)|}
\leq
\left( {3\over 4} \right) ^{\sqrt{\lg (n)} },
$$
and the right hand side goes to zero as $n\to \infty $.

Therefore, the Automorphism Conjecture is true for this
class of ordered sets.
\qed

\begin{lem}
\label{wid12AClem}

There is an $N\in {\mat N} $ such that, for every
ordered set $P$
of width $w\leq 11$ with
$n\geq N$ elements,
no nontrivial order-autonomous antichains,
and at most $\lg (n) $
elements
contained in
max-locked interdependent orbit unions of height $1$,
we have that
${|{\rm Aut} (P)|
\over
|{\rm End} (P)|}
\leq 2^{-0.005 n} $.

\end{lem}

{\bf Proof.}
By
Lemma \ref{orbits11orless},
for all tight interdependent orbit unions $(U,{\cal D} )$
such that, for all $D\in {\cal D}$ we have $|D|\leq 11$, and all
orbits $|D|$ such that $\Lambda _{\cal D} (D)$ contains $A_{|D|} (D)$
satisfy $|D|\leq 5$, we have that
$|{\rm Aut} _{{\cal D} } (U )|\leq 2^{0.993|U |} $.

Let $(U, \cal D)$ be a tight
interdependent orbit union such that there is an
orbit $D$ with more than $5$ elements such that
$|\Lambda _{\cal D} (D)|$ contains $A_{|D|} (D)$.
By Lemma \ref {orbits11orlessmaxlocked},
$(U, \cal D)$ is
max-locked of width $v\leq w\leq 11$.
Now either $U$ has height $1$, or,
a quick computation with factorials shows that
$|{\rm Aut} _{\cal D} (U)|=
v!\leq 2^{2.4v}\leq 2^{0.8|U|} $.

Let
$Z\subseteq P$ be the set of all points in $P$
that are
contained in (possibly trivial) interdependent orbit
unions that are
not max-locked of height $1$.
Let $K:=|Z|$.
By assumption $n-K\leq \lg (n) $.
By the above,
the number of
restrictions of automorphisms of $P$
to
$Z$ is
bounded by
$2^{0.993|Z |} $.
The number of
restrictions of automorphisms of $P$
to $P\setminus Z$
is bounded by
$|P\setminus Z|! = (n-K)!\leq \lg (n)!\leq
\lg (n)^{\lg (n)} =2^{\lg (n) \lg (\lg (n))}$.

Hence,
the number of automorphisms of
$P$ is bounded by
$
2^{0.993n}
2^{\lg (n) \lg (\lg (n))}
$,
which, for large enough $n$, is smaller than
$
2^{0.994n}
$.
The claim now follows from
Lemma \ref{2tonforbddwidth}
with $\varepsilon :=0.001$.
\qed

\vspace{.1in}

We are left with the task to consider
ordered sets with nontrivial order-autonomous antichains.
Proposition \ref{onelexiter} below focuses on the
more general lexicographic sum construction.

\begin{define}
\label{lexsumdef}

(See, for example, \cite{HH1,MR}.)
Let $T$ be a nonempty ordered set considered as an index set.
Let
$\{ P_t \} _{t\in T} $ be a family of
pairwise
disjoint nonempty ordered sets that
are all disjoint from $T$. We define
the {\bf lexicographic sum} $L\{ P_t \mid t\in T\} $
{\bf (of the $P_t $ over $T$)}
to be
the union
$\bigcup _{t\in T} P_t $
ordered by letting $p_1 \leq p_2 $ iff
either
\begin{enumerate}
\item
There are distinct
$t_1 , t_2 \in T$ with
$t_1 <t_2 $,
such that
$p_i\in P_{t_i } $,
or
\item
There is a $t\in T$ such that
$p_1 , p_2 \in P_t $ and $p_1 \leq  p_2 $ in $P_t $.
\end{enumerate}
The ordered sets
$P_t $ are the {\bf pieces} of the lexicographic sum and
$T$ is the {\bf index set}.

\end{define}

\begin{prop}
\label{onelexiter}

Let ${\cal C}$ be a class of finite
ordered sets for which the Automorphism
Conjecture holds and let
${\cal L} _{\cal C}$ be the
class of lexicographic sums $L\{ P_t |t\in T\} $ such that all $P_t \in {\cal C}$
and
$T\in {\cal C} $ is indecomposable or a chain, or,
all $P_t $ are antichains
and $T\in {\cal C} $ does not contain any nontrivial order-autonomous antichains.
Then
the Automorphism
Conjecture holds for
${\cal L} _{\cal C}$.

\end{prop}

{\bf Proof.}
For every $n\in {\mat N}$, let
$$h(n):=
\max _{P\in {\cal C} , |P|\geq n}
{ |{\rm Aut} (P)|\over |{\rm End} (P)| } .
$$
Because
$\lim _{n\to \infty } h(n)=0$,
for $n\geq 2$, we have that
$h(n)<1$.

For every
$L\{ P_t |t\in T\} \in {\cal L} _{\cal C} $,
let
$A(P_t:t\in T):=
\sum _{|P_t|>1} |P_t |$.

Let
$\varepsilon >0$.
Fix
$K\in {\mat N}$ such that,
for all $k\geq K$, we have that
$(h(2))^k
<\varepsilon $ and
$h(k)
<\varepsilon $.
Choose $N\in {\mat N}$ such that
$N\geq K^2 $ and such that, for all
$n\geq N$, we have
$h\left( n-
K^2 \right) \leq \varepsilon ^2 $
and $\left( K^2 \right) !\leq {1\over \sqrt{h\left( n-
K^2 \right) } }$.
Then, for every
$L\{ P_t |t\in T\} \in {\cal L} _{\cal C} $
with $n\geq N$ elements, we obtain the following.

Recall that automorphisms map
order-autonomous sets to order-autonomous sets.
Therefore,
for every automorphism $\Phi \in {\rm Aut} (P)$,
all endomorphisms of all sets $P_s $ translate into
pairwise distinct
order-preserving maps from $P_s $ to $\Phi [P_s ]$.
In case
$A(P_t:t\in T) >K^2 $,
this provides the following estimates.

In case
$A(P_t:t\in T) >K^2 $ and
there is an $s\in T$ with
$|P_s |>K$, we have
\begin{eqnarray*}
{
|{\rm Aut} (L\{ P_t |t\in T\} )|
\over
|{\rm End} (L\{ P_t |t\in T\} )|
}
& \leq &
{
|{\rm Aut} (P_s )|
\over
|{\rm End} (P_s )|
}
\leq h(|P_s | )
<\varepsilon
\end{eqnarray*}

In case
$A(P_t:t\in T) >K^2 $ and, for all
$t\in T$ we have
$|P_t |\leq K$,
let $R$ be the set of
indices $r\in T$ such that
$|P_r |>1$.
Then $|R|\geq K$
and we have
\begin{eqnarray*}
{
|{\rm Aut} (L\{ P_t |t\in T\} )|
\over
|{\rm End} (L\{ P_t |t\in T\} )|
}
& \leq &
\prod _{r\in R}
{
|{\rm Aut} (P_r )|
\over
|{\rm End} (P_r )|
}
\\
& \leq &
(h(2))^K
<\varepsilon
\end{eqnarray*}

This leaves the
case
$A(P_t:t\in T) \leq K^2 $.
Because
$T\in {\cal C} $ is indecomposable or a chain or
$T\in {\cal C} $ does not contain any nontrivial order-autonomous antichains and
all $P_t $ are antichains,
automorphisms map sets $P_t $ to sets $P_t $ and every automorphism
induces a corresponding automorphism on $T$.
Hence we obtain the following.
\begin{eqnarray*}
{
|{\rm Aut} (L\{ P_t |t\in T\} )|
\over
|{\rm End} (L\{ P_t |t\in T\} )|
}
& \leq &
{
|{\rm Aut} (T)|\cdot A(P_t:t\in T)!
\over
|{\rm End} (T )|
}
\\
& \leq &
h\left(
|T|
\right)
A(P_t:t\in T)!
\\
& \leq &
h\left(
|T|
\right)
\left( K^2 \right) !
\\
& \leq &
\sqrt{h\left(
n
-
K^2
\right) }
<\varepsilon
\end{eqnarray*}
\qed

\begin{theorem}
\label{ACw<=12}

The
Automorphism Conjecture is true for ordered sets of
width $w\leq 11$.

\end{theorem}

{\bf Proof.}
By
Proposition \ref{manyinmaxlock}
and Lemma \ref{wid12AClem}, the
Automorphism Conjecture is true for
ordered sets of
width $w\leq 11$ without nontrivial order-autonomous antichains.
Now apply Proposition \ref{onelexiter},
using the
ordered sets of width $\leq 11$ without order-autonomous antichains
and the antichains as the base class ${\cal C}$.
This proves the Automorphism Conjecture for a class
that contains all
ordered sets of
width $w\leq 11$.
\qed

\begin{remark}

{\rm
By Table \ref{badgrouptable}, the only primitive group
of degree $|X|=13$ satisfies $|G|\leq 2^{0.9582|X|} $.
Thus orbits of size 13 could be added to the key
Lemma \ref{orbits11orless}.
However, for the imprimitive groups of
degree $|X|=12$, the order
$(6!)^2 2!$
exceeds $2^{1.5|X|} $.
It should not be too hard to
obtain an analogue of Lemma \ref{orbits11orlessmaxlocked}
to manage this permutation group and to make sure that there
are no other problematic cases.
However, to establish the Automorphism Conjecture for all
ordered sets of width $\leq 13$, we would also need to
treat the situation of orbits
as in Example \ref{noDendos} below (with $M=6$).
Hence such an attempt at an extension will
encounter additional challenges.

Although this possible extension
feels very feasible
(also see Remark \ref{gradedto13} below), and the author is not superstitious,
it stands to reason that energy expended on small orbits
should be expended as indicated in Remark
\ref{useQinfuture}:
Strengthenings of Lemma \ref{expolem} are likely to have
more far-reaching consequences, as
can be seen from Section \ref{primnestsec}.
}

\end{remark}

\begin{remark}
\label{gradedto13}

{\rm
An ordered set $P$ is called {\bf graded}
iff there is a function $g:P\to {\mat N}$ such that,
if $y$ is an upper cover of $x$, then
$g(y)-g(x)=1$.
For any automorphism $\Phi :P\to P$ of a graded ordered set
$P$, and for all $x\in P$, we have that
$g(\Phi (x))=g(x)$.
Hence, for a graded ordered set $P$,
Lemma \ref{bipartiteOG} can be applied
with
$B:=\{ x\in P: g(x)=2k, k\in {\mat N}\} $
and
$T:=\{ x\in P: g(x)=2k-1, k\in {\mat N}\} $.
Moreover, the computation in
Lemma \ref{controllednest} can be done for primitive nestings in which no
factor
with ${|B_{j+1} |\over |B_j |}\geq 7$
contains an alternating group with
${|B_{j+1} |\over |B_j |}$ elements to obtain a bound of
$|G |\leq 2^{1.95|X|} $.
Hence, for a {\em graded} tight
interdependent orbit union
$(U,{\cal D} )$, we can
obtain a version of Theorem \ref{permgrtoordaut}
that guarantees
$
|{\rm Aut} _{\cal D} (U) |
\leq 2^{{1\over 2}1.95|U |} $
under the condition that no factor
with ${|B_{j+1} |\over |B_j |}\geq 7$
in a primitive nesting of a
$\Lambda _{\cal D} (D )$
contains an alternating group with
${|B_{j+1} |\over |B_j |}$ elements.
Now the same argument as given in this section proves the
Automorphism Conjecture for {\em graded} ordered sets
of width at most 13.
}

\end{remark}

\section{More Endomorphisms for Heights 2 and 3}
\label{moreendosec}

Theorem \ref{ACw<=12} clearly demonstrates the utility of the trivial
realization that having more endomorphisms is better when considering the
Automorphism Conjecture:
The guarantee of
at least
$2^{\left( 1
-\varepsilon \right) n} $
endomorphisms in Lemma \ref{2tonforbddwidth} allows
for slightly less stringent estimates for the number
of automorphisms to carry the day.
Remark \ref{smallprimitivefactors} indicates that
exclusion of a few more factors in primitive nestings even
guarantees $|{\rm Aut} _{{\cal D} } (U )|\leq
2^{0.76|U|} .$
In Corollary \ref{primbounds} below, we will encounter a
natural
situation in which
$|{\rm Aut} _{{\cal D} } (U )|\leq
2^{0.69|U|} .$

Theorem 1 of \cite{DRSW} guarantees that
every ordered set $P$
has at least $2^{{2\over 3} |P|} $ endomorphisms.
Clearly, this is a ``near miss" when
$|{\rm Aut} _{{\cal D} } (U )|\leq
2^{0.69|U|} .$
Theorem 1 of \cite{DRSW} is proved by showing that
every ordered set $P$ of height $1$ has at least
$2^{|P|} $ endomorphisms (see (C) in Section 2 of \cite{DRSW}),
and
that every ordered set of height $h$ has at least
$2^{{h\over h+1} n} $ endomorphisms (see (D) in Section 2 of \cite{DRSW}).
Hence we can improve the lower bound on the number of endomorphisms
by improving the lower bounds for heights
$2$ and $3$.
The improvement for height $2$
suffices for the situation in Corollary \ref{primbounds}.
Since the
factor in the
improvement exceeds the factor
${3\over 4} $ for height $3$, it is natural to
go one step further to height $3$.

\begin{lem}
\label{endht2}

Let $P$ be an ordered set of height $2$ with $n$ elements.
Then we have
${ |{\rm Aut} (P)|\over |{\rm End} (P)| }
\leq 2^{\sqrt{{1\over 2} \lg (n) } } $,
or, for $n\geq 22$, we have
$|{\rm End} (P)|\geq 2^{{\lg (3)\over 2} n} >2^{0.7924n} $.

\end{lem}

{\bf Proof.}
Let
$n_j $ denote the number of
elements of rank $j$.
If $P$ has more than
$\lg (n)$ irreducible points, then, by
Theorem 1 in \cite{LiuWan},
we have
${ |{\rm Aut} (P)|\over |{\rm End} (P)| }
\leq 2^{\sqrt{{1\over 2} \lg (n) } } $.
We are left with the case that
$P$ has fewer than
$\lg (n)$ irreducible points.

If all non-maximal
elements of rank $1$ are irreducible, then
$n_1 <\lg (n)$ and $P\setminus R_1 $
is an ordered set of height $1$
with more than
$n-\lg(n) $ elements and hence more than
$2^{n-\lg (n)} $ endomorphisms.
For $n\geq 22$, we have $2^{n-\lg (n)} \geq 2^{{\lg (3)\over 2}n} $.
We are left with the case that
there is a non-irreducible non-maximal
element $c$ of rank $1$.

Let $b_1 ,b_2 <c$ be two lower covers of $c$
and let
$t_1 ,t_2 >c$
be two upper covers of $c$.
Mapping all elements of $R_1 $ to $c$,
mapping the elements of $R_0 $
in any fashion to $\{ b_1 , b_2 ,c \} $,
and mapping the elements of $R_2 $
in any fashion to $\{ c,t_1 , t_2 \} $
generates
$3^{n_0 +n_1 } =2^{\lg (3)(n_0 +n_2 )} $
endomorphisms.
Hence, in case
$n_0 +n_2 \geq
{n\over 2} $,
we have at least
$2^{{\lg (3)\over 2}  n} $
endomorphisms.
We are left with the case that
$n_0 +n_2 <
{n\over 2} $.

In this case, we have
$n_1 \geq
{n\over 2} $.
Mapping
$R_0 $ to $b_1 $, $R_2 $ to $t_1 $ and
mapping the elements of $R_2 $ to
$\{ b_1 ,c, t_1 \} $ in any fashion
generates
$3^{n_1 } \geq 2^{{\lg (3)\over 2}  n} $
endomorphisms.
\qed

\begin{lem}
\label{endht3}

Let $P$ be an ordered set of height $3$ with $n$ elements.
Then we have
${ |{\rm Aut} (P)|\over |{\rm End} (P)| }
\leq 2^{\sqrt{{1\over 2} \lg (n) } } $,
or, for $n\geq 27$,
we have $|{\rm End} (P)|\geq 2^{{\lg (3)\over 2} (n-\lg(n))}$,
which, for $n\geq 180,000$, exceeds
$2^{0.7924n} $.

\end{lem}

{\bf Proof.}
Let
$n_j $ denote the number of
elements of rank $j$.
If $P$ has more than
$\lg (n)$ irreducible points, then, by
Theorem 1 in \cite{LiuWan},
we have
${ |{\rm Aut} (P)|\over |{\rm End} (P)| }
\leq 2^{\sqrt{{1\over 2} \lg (n) } } $.
We are left with the case that
$P$ has fewer than
$\lg (n)$ irreducible points.

Consider the case that
all non-maximal
elements of rank $2$ are irreducible.
In this case,
$n_2 <\lg (n)$ and $P\setminus R_2 $
is an ordered set of height $2$
with more than
$n-\lg(n) \geq 22 $ elements.
The image of an irreducible point under an automorphism is determined by
its unique upper or lower cover.
Therefore $|{\rm Aut} (P)|=|{\rm Aut} (P\setminus R_2 )|$.
Now, via Lemma \ref{endht2}, we obtain
${ |{\rm Aut} (P)|\over |{\rm End} (P)| }
\leq
{|{\rm Aut} (P\setminus R_2 )|\over |{\rm End} (P\setminus R_2 )| }
\leq 2^{\sqrt{{1\over 2} \lg (|P\setminus R_2 |) } }
\leq 2^{\sqrt{{1\over 2} \lg (n) } }
$,
or, for $n\geq 27$,
$|{\rm End} (P)|\geq |{\rm End} (P\setminus R_2 )|
\geq
2^{{\lg (3)\over 2} (|P\setminus R_2 |)}
\geq
2^{{\lg (3)\over 2} (n-\lg(n))}$.
Similarly, if all non-minimal
elements of dual rank $2$ are irreducible, then
there are fewer than $\lg (n)$
elements of dual rank $2$,
$P\setminus
\{ x\in P\setminus \min(P): {\rm dual \ rank} (x)=2\}  $
is an ordered set of height $2$
with more than
$n-\lg(n) $ elements, and, again,
we obtain
the
requisite bounds via Lemma \ref{endht2}.
We are left with the case that
there is a
non-irreducible non-maximal
element $c$ of rank $2$
and that there is a
non-irreducible non-minimal
element $d$ of dual rank $2$.

Let $t_1 , t_2 >c$ be two upper covers of $c$
and let $a<b<c$.

If $n_2 +n_3 \geq {n\over 2}$,
mapping
$R_0 \cup R_1 $ to $a$, mapping
$R_2 $ to $\{ a,b,c\} $, and
mapping
$R_3 $ to $\{ c,t_1,t_2 \} $
generates $3^{n_2 +n_3 } \geq 2^{ {\lg (3)\over 2} n} $
endomorphisms.
If $n_1 +n_3 \geq {n\over 2}$,
mapping
$R_0 $ to $a$, mapping
$R_2 $ to $c$, mapping
$R_1 $
to $\{ a,b,c\} $, and
mapping
$R_3 $ to $\{ c,t_1,t_2 \} $
generates $3^{n_1 +n_3 } \geq 2^{ {\lg (3)\over 2} n} $
endomorphisms.
If $n_0 +n_3 \geq {n\over 2}$,
mapping
$R_1 \cup R_2 $ to $c$, mapping
$R_0 $ to $\{ a,b,c\} $, and
mapping
$R_3 $ to $\{ c,t_1,t_2 \} $
generates $3^{n_0 +n_3 } \geq 2^{ {\lg (3)\over 2} n} $
endomorphisms.

If there is an
$i\in \{ 1,2,3\} $ such that
$n_0 +n_i \geq {n\over 2} $, then,
dual to the preceding argument, we can
use $d$, two of its lower covers, and a chain $d<b<a$,
to construct
$\geq 2^{ {\lg (3)\over 2} n} $
endomorphisms.

Using
the larger of $n_0 $ and $n_3 $ and the larger of
$n_1 $ and $n_2 $, we conclude that,
for some $j\in \{ 0,1,2\} $, we must have
$n_j +n_3 \geq {n\over 2} $, or,
for some $i\in \{ 1,2,3\} $, we must have
$n_0 +n_i \geq {n\over 2} $.
This, and a quick computation to verify the
final claim,
concludes the proof.
\qed

\begin{theorem}
\label{moreendos}

Let $P$ be an ordered set
with $n$ elements.
Then we have
${ |{\rm Aut} (P)|\over |{\rm End} (P)| }
\leq 2^{\sqrt{{1\over 2} \lg (n) } } $,
or, for $n\geq 27$,
we have $|{\rm End} (P)|\geq 2^{{\lg (3)\over 2} (n-\lg(n) )}$,
which, for $n\geq 180,000$, exceeds
$2^{0.7924n} $.

\end{theorem}

{\bf Proof.}
This follows from
(C) in Section 2 of \cite{DRSW}, which says that
every ordered set of height $1$ with $n$ elements has at least
$2^n $ endomorphisms,
from
Lemmas \ref{endht2} and \ref{endht3}
for ordered sets of heights $2$ and $3$, respectively, and
from (D) in Section 2 of \cite{DRSW},
which says that every ordered set of height $h(\geq 4)$
with $n$ elements has at least
$2^{{h\over h+1} n} $ endomorphisms.
\qed

\section{Orbits Carrying Primitive Permutation Groups}
\label{primorbsec}

Section \ref{primnestsec} shows the significant utility of the
ideas presented here beyond the proof of concept in
Section \ref{ACsmwidthproof}.
The only drawback to
Section \ref{primnestsec} is that there is no easily stated
class of ordered sets for which the Automorphism Conjecture is confirmed.
In this section we will focus on
tight interdependent orbit unions such that
all orbits carry
primitive permutation groups.
We will see that the Automorphism Conjecture can be
confirmed for ordered sets whose
natural interdependent orbit unions are so that all
orbits carry primitive permutation groups,
and which do not exhibit the problem indicated in
Example \ref{noDendos}.
In particular, this gives a first insight how to possibly handle
factorial factors.

\begin{prop}
\label{primshouldbeok}

Let $(U,{\cal D} )$ be a tight interdependent orbit union
such that all $\Lambda _{\cal D} (D)$ are primitive.
Then, for any $D\in {\cal D}$, we have
$|{\rm Aut} _{{\cal D} } (U )|= |\Lambda _{\cal D} (D) |$.

\end{prop}

{\bf Proof.}
We will prove by induction, that, for any
subgroup
$G^* $ of $|{\rm Aut } _{{\cal D}} (U)|$ such that
all $\Lambda _{\cal D} ^* (D)$ are primitive,
for any $D\in {\cal D}$, we have that,
$|G^* |=|\Lambda ^* (D)|$.

The proof is an induction on the number $m=|{\cal D}|$ of
orbits.
For the base case $m=2$, first note that, because
order-autonomous antichains in $D$ constitute blocks of
$\Lambda _{\cal D} ^* (D)$,
no orbit $D\in {\cal D}$ can contain a
nontrivial order-autonomous antichain.
In the absence of
nontrivial order-autonomous antichains within either of
$D_1 $ and $D_2 $,
by Lemma \ref{allbutoneac},
every $\Phi \in {\rm Aut} _{\cal D} (D)$
is completely determined by $\Phi |_{D_1 } $ or by $\Phi |_{D_2 } $,
which leads to
$|G^* |=|\Lambda ^* (D_i )|$
for either $i\in \{ 1,2\} $.

For the induction step $\{1, \ldots , m-1\} \to m$,
let ${\cal D} \in {\cal D} $.
We first note that ${\cal O} (U,{\cal D})$ has at least two
noncutvertices.
Therefore, we can relabel the
orbits
in
${\cal D}$ such that $D_m $ is a noncutvertex that is not equal to $D$.
By Theorem \ref{pruneorbit6},
we have
$
|G^*|
=
\left| {\rm Aut } _{{\cal D}_n} ^* (U_n )\right|
\left|{\rm Aut } _{{\cal D}_Q } ^* (Q)\right| $.
Because all
$\Lambda _{\cal D} ^* (D)$ are primitive,
all partitions ${\cal A} ^j $ from
Notation \ref{blocksinDj} must be partitions into singletons.
Therefore,
the only element of
${\rm Aut} _{{\cal D} _Q } ^* (Q)$
is the identity.
Moreover, for all
$j\in \{ 1, \ldots , m-1\} $, we have that
the permutation group induced by
${\rm Aut } _{{\cal D}_m} ^* (U_m )$ on $D_j $ is
$\Lambda _{\cal D} ^* (D_j )$ and hence it is primitive.
Therefore we can apply the induction hypothesis to
${\rm Aut} _{{\cal D} _Q } ^* (Q)$.

We conclude that
$
|G^* |
=
\left| {\rm Aut } _{{\cal D}_m} ^* (U_m )\right|
=
\left| \Lambda ^* (D)\right|
$, which concludes the induction.
The result now follows by using $G^* ={\rm Aut} _{\cal D} (U)$
in the claim we just proved.
\qed

\begin{cor}
\label{primbounds}

Let $(U,{\cal D} )$ be a tight interdependent orbit union
such that all $\Lambda _{\cal D} (D)$ are primitive and don't contain
$A_{|D|} (D)$.
Then
$|{\rm Aut} _{{\cal D} } (U )|\leq
2^{1.38\min \{ |D|:D\in {\cal D} \} } \leq 2^{0.69|U|} .$

\end{cor}

{\bf Proof.}
By Theorem \ref{marotiexpbound}, we obtain that, for all $D\in {\cal D}$,
we have
$|\Lambda _{\cal D} (D)| < 2^{1.38|D|}$.
By Proposition \ref{primshouldbeok},
we obtain that, for any $D\in {\cal D}$, we have
$|{\rm Aut} _{{\cal D} } (U )|= |\Lambda _{\cal D} (D) |< 2^{1.38|D|}$.
Because there is a $D\in {\cal D}$ with $|D|\leq {1\over 2} |U|$, we
obtain the conclusion.
\qed

\begin{prop}
\label{auffaecherung}

Let $(U,{\cal D} )$ be a tight interdependent orbit union
such that all $\Lambda _{\cal D} (D)$ are primitive,
such that there is a $D\in {\cal D}$ such that
$\Lambda _{\cal D} (D)$ contains a copy of $A_{|D|} (D)$ and such that not
all
orbits
$C\in {\cal D}$ have the same size.
Then $|U|\geq
{1\over 2} |D|(|D|+1)$ and hence
$|{\rm Aut} _{\cal D} (U)|\leq 2^{{1\over 2} |U|} $.

\end{prop}

{\bf Proof.}
Let $D$ be an
orbit
such that
$\Lambda _{\cal D} (D)$ contains a copy of $A_{|D|} (D)$.
By Lemma \ref{firstauffaecherung} (and trivially for $|D|=2$),
if $D\updownharpoons _{\cal D} E$
and $|D|=|E|=:w$, then
$D\cup E$ is isomorphic to an ordered set
$wC_2 $ or $S_w $ and hence
$\Lambda _{\cal D} (E)$ contains a copy of $A_{|E|} (E)$.
Because not all
frames in
${\cal D}$
are of the same size, we can assume,
without loss of
generality, that there is an
$E\in {\cal D}$ such that
$D\updownharpoons _{\cal D} E$
and $|D|\not= |E|$.
Because $\Lambda _{\cal D} (E)$ is primitive,
we obtain
$|D|\geq 3$.
By Lemma \ref{firstauffaecherung},
we have that there is a $k\in \{ 2, \ldots , |D|-2\} $
such that
$|E|=\pmatrix{ |D|\cr k\cr }
\geq {1\over 2} |D|(|D|-1)$.

Hence
$|U|\geq |D|+{1\over 2} |D|(|D|-1)={1\over 2} |D|(|D|+1)$.
By Proposition \ref{primshouldbeok}, we have
$|{\rm Aut} _{{\cal D} } (U )|\leq |\Lambda _{\cal D} (D)|
\leq |D|!\leq 2^{0.4605\cdot {1\over 2} |D|(|D|+1) }
\leq 2^{{1\over 2} |U|} $.
\qed

\begin{remark}

{\rm
Theorem \ref{moreendos},
Corollary \ref{primbounds} and
Proposition \ref{auffaecherung}
show that, when considering the Automorphism Conjecture,
for the case of
$(U,{\cal D} )$ being a tight interdependent orbit union
such that all $\Lambda _{\cal D} (D)$ are primitive,
the only case that remains is the case that
there is a $D$ such that
$\Lambda _{\cal D} (D)$ contains a copy of $A_{|D|} (D)$ and such that
all
orbits
$C\in {\cal D}$ have the same size.

More progress is possible
with the ideas from Lemma \ref{maxlockh1presrank}
and Proposition \ref{manyinmaxlock}:
When the orbit
graph is a tree, it can be proved that
$|{\rm End}_{\cal D} (U)|\geq (w-1)^w $,
where $w$ is the width of an orbit.
However, although
Lemma \ref{self-interference} below
shows that the structure
must be very specific
when the orbit graph contains cycles,
Example \ref{noDendos} shows that ${\cal D}$-endomorphisms
will not provide the complete answer here.
Although this may be a bit anticlimactic,
Example \ref{noDendos} is representative for the only
remaining problem for these interdependent orbit unions, which
should be resolvable in the future.
}

\end{remark}

\begin{define}

Let $(U,{\cal D})$ be an interdependent orbit union
and let
${\cal C} =\{
D_1 \updownharpoons _{\cal D} D_2 \updownharpoons _{\cal D} \cdots
\updownharpoons _{\cal D} D_N \updownharpoons _{\cal D} D_1 \} $
be a cycle in ${\cal O} (U,{\cal D} )$.
Then ${\cal C}$ is called a {\bf lock cycle}
iff there is an $M>1$ such that,
for all $i\in \{ 1, \ldots , N\} $
and with index arithmetic modulo $N$,
we have that
$D_i \cup D_{i+1} $ is isomorphic to
$S_M $ or to $MC_2 $.
For
$x,y\in D_1 $, we write
$x\circlearrowright y$
and say $y$ is {\bf cycle locked} to $x$ iff there is a sequence
$x=x_1 , x_2 , \ldots , x_N , x_{N+1} =y$ such that,
for all $i\in \{ 1, \ldots , N\} $, we have that $x_i \in D_i $, and,
if $D_i \cup D_{i+1} $ is isomorphic to
$S_M $, then $x_i \not\sim x_{i+1} $, and,
if $D_i \cup D_{i+1} $ is isomorphic to
$MC_2 $, then $x_i \sim x_{i+1} $.

\end{define}

\begin{lem}
\label{self-interference}

Let $(U,{\cal D})$ be an interdependent orbit union.
Let
$D_1 \updownharpoons _{\cal D} D_2
\linebreak
\updownharpoons _{\cal D} \cdots
\updownharpoons _{\cal D} D_N \updownharpoons _{\cal D} D_1 $
be a lock cycle in ${\cal O} (U,{\cal D} )$
such that
there are distinct $x,y\in D_1 $ such that
$x\circlearrowright y$.
Then
$|\Lambda _{\cal D} (D_1 )|\leq {1\over M-1} M!$.
In particular, $\Lambda _{\cal D} (D_1 )$ does not contain a copy of
$A_{|D_1 |} (D_1 )$.

\end{lem}

{\bf Proof.}
Every $\Phi \in {\rm Aut} _{\cal D} (U)$ must preserve the
$\circlearrowright $-relation and every $x\in D_1 $ is
$\circlearrowright $-related
to exactly one element in $D_1 $. Therefore, if
$x\circlearrowright y$ and $x\not= y$, then,
for
every $\Phi \in {\rm Aut} _{\cal D} (U)$,
$\Phi (x)$ uniquely determines $\Phi (y)$ as the unique
element that is
$\circlearrowright $-related
to $\Phi (x)$.
Hence $|\Lambda _{\cal D} (D_1 )|\leq {1\over M-1} M!$.
\qed

\vspace{.1in}

Lemma \ref{self-interference} gives further information about the
structure of interdependent orbit unions with
induced permutation groups
$\Lambda _{\cal D} (D)$ that contain
the alternating group on $D$ and such that all orbits are of the same size:
In every lock cycle, all elements are cycle locked to themselves.
Unfortunately, we cannot easily dispense with this
remaining case in the same way we did
when the width was bounded by 11 because of the following
situation (and similar ones).

\begin{exam}
\label{noDendos}

{\rm
{\em
Let $M\in {\mat N}\setminus \{ 1,2\} $ and consider the
interdependent orbit union
$(U, {\cal D})$ with ${\cal D} =\{ D_1 , D_2 , D_3 , D_4 \} $,
$D_i =\{ x_1 ^i , \ldots , x_M ^i \} $ such that
$x_j ^1 <x_k ^2 $ iff
$j\not= k$,
$x_j ^3 <x_k ^2 $ iff
$j= k$,
$x_j ^3 <x_k ^4 $ iff
$j= k$,
$x_j ^1 <x_k ^4 $ iff
$j= k$, and no further comparabilities.
Then every $f\in {\rm End} _{\cal D} (U)$
is an automorphism.
}

Let
$f\in {\rm End} _{\cal D} (U)$.
Without loss of generality, we can assume that
$f$ is a retraction onto its range.
Hence $f|_{D_1 \cup D_2 } $ is a retraction that preserves minimal and maximal elements.
Let $x_j ^1 \in f[D_1 ]$. Then no element of
$D_2 \setminus \{ x_j ^2 \} $ is mapped to
$x_j ^2 $.
Similarly, if
$x_j ^2 \in f[D_2 ]$, then no element of
$D_1 \setminus \{ x_j ^1 \} $ is mapped to
$x_j ^1 $.

Let $x_j ^1 \in f[D_1 ]$.
Then $x_j ^4 \in f[D_4 ]$,
$x_j ^3 \in f[D_3 ]$
and
$x_j ^2 \in f[D_2 ]$.
Similarly,
if $x_j ^2 \in f[D_2 ]$,
then $x_j ^1 \in f[D_1 ]$.
Hence
$x_j ^1 \in f[D_1 ]$
iff
$x_j ^2 \in f[D_2 ]$.
Via the preceding paragraph,
we conclude that no element in $D_1 \setminus f[D_1 ]$ can be mapped to any
element of $f[D_1 ]$.
Hence $f[D_1 ]=D_1 $ and $f$ must be an automorphism.
\qex
}

\end{exam}

\section{Conclusion}
\label{conclusec}

The insights presented in Sections \ref{orbclustsec} and \ref{orbgraphsec}
are a significant addition to our
understanding of the automorphic structure of
ordered sets.
The connection to permutation groups on orbits in
Section \ref{AuttoLambda} is a significant advance
towards confirming the Automorphism Conjecture.
The remarks throughout this presentation indicate a
wide variety of fruitful paths forward.
Remark \ref{useQinfuture} indicates that
deeper use of the
order-theoretical properties of the
residual
sets $(Q,{\cal D} _Q )$
should allow a further
tightening of the bounds so far.
Remark \ref{smallprimitivefactors} articulates that the main remaining
challenges will be with factor groups of the
$\Lambda _{\cal D} (D)$
which contain a subgroup $A_{|D|} (D)$
and possibly with some factor groups from
among the
small primitive groups from
Table \ref{badgrouptable}.
Remark \ref{refinekeylemmss}
indicates that further refinements of the key
Lemma \ref{expolem}
are a worthy target for future investigations.
The work in Section \ref{primorbsec} shows that further progress can be made when
the problematic factors in
$\Lambda _{\cal D} (D)$
necessitate that certain orbits $E\updownharpoons _{\cal D} D$
must be large. Although Section \ref{primorbsec}
focuses on primitive groups $\Lambda _{\cal D} (D)$,
a group $G\hochdamit ^{G\cdot B} _A $, where $A,B$ are blocks
of $\Lambda _{\cal D} (D)$, that contains
an alternating group on the blocks in
$A[G\cdot B]$
will, through a variation on
Lemma \ref{firstauffaecherung},
for every $E\updownharpoons _{\cal D} D$,
require that $D\cup E$
contains two levels of a Boolean algebra or a set $wC_2 $, but
these levels will be made up
of blocks rather than points, and there may be multiple parallel copies.
Unfortunately, this does not automatically create the
size imbalance
that was used in
Proposition \ref{auffaecherung}, because the blocks in $E$ might contain
structures that, in turn, require
larger blocks in another level of
the primitive nesting of $D$.
This and the fact that, aside from the alternating and symmetric groups,
all remaining problematic groups are small, and can indeed
be explicitly identified in a short list, makes the
author very hopeful for the
future of this work.

\vspace{.1in}

{\bf Acknowledgement.}
The author very much thanks Frank a Campo for a thorough
reading of an earlier, much more technical, version of this paper,
for the detection of
mistakes, many helpful suggestions,
improvements of some proofs, and
for the suggestion to
more accurately call the elements of ${\cal D}$ ``frames" instead of ``dictated orbits."

\section{Declarations}

\begin{itemize}
\item
No funds, grants, or other support was received.

\item
The author has no relevant financial or non-financial
interests to disclose.

\item
Data sharing is not applicable, as no datasets
were generated for this paper.
The simple R code mentioned in the proof of Lemma \ref{controllednest}
can be made available upon request.

\end{itemize}


\begin{thebibliography}{99}

\bibart{BMgroupAutP}
{J. A. Barmak, E. G. Minian}{2009}
{Automorphism Groups Of Finite Posets}
{Discrete Mathematics}{309}{3424--3426}

\bibitem{BonaMartin}
M. B\'ona
and R. Martin (2022),
The endomorphism conjecture for graded
posets of width 4,
\url{https://arxiv.org/abs/2205.15378}


\bibart{BuekLee}
{F. Buekenhout, D. Leemans}{1996}
{On the List of Finite Primitive Permutation Groups
of Degree $\leq 50$}
{J. Symbolic Computation}{22}{215--225}


%
%





\bibbook{Dixmort}
{J. D. Dixon, B. Mortimer}{1996}
{Permutation Groups}
{Springer}{New York, NY}






\bibart{DRSW}
{D. Duffus, V. R\"odl, B. Sands, R. Woodrow}{1992}
{Enumeration of order-preserving maps}{Order}{9}{15-29}







\bibart{HH1}
{H. H\"oft and M. H\"oft}{1976}
{Some fixed point
theorems for partially ordered sets}{Can. J. Math.}{28}{992--997}









\bibart{LRZ}
{W.-P. Liu, I. Rival and N. Zaguia}{1995}
{Automorphisms, Isotone self-maps and cycle-free orders}
{Discrete Mathematics}{144}{59--66}

\bibart{LiuWan}
{W.-P. Liu and H. Wan}{1993}
{Automorphisms and Isotone Self-Maps of Ordered Sets with
Top and Bottom}{Order}{10}{105--110}

\bibart{MR}
{M. Malicki and A. Rutkowski}{2004}
{On operations and linear extensions of well partially ordered sets}
{Order}{21}{7--17}



\bibart{Maroti}
{A. Mar\'oti}{2002}
{On the orders of primitive groups}
{Journal of Algebra}{258}{631--640}






\bibart{Proe}
{H. J. Pr\"omel}{1987}
{Counting unlabeled structures}
{J. Comb. Theory Ser. A}{44}{83--93}




\bibcoll{RiRut}
{I. Rival and A. Rutkowski}{1991}
{Does almost every isotone self-map have a
fixed point?}{Bolyai Math. Soc.}{Extremal Problems for Finite
Sets}{Bolyai Soc. Math. Studies 3,
Vis\'egrad, Hungary}{413-422}






\bibbook{Schbook}
{B. Schr\"oder}{2016}{Ordered Sets -- An Introduction with Connections
from Combinatorics to Topology (second edition)}
{Birkh\"auser Verlag}{Boston, Basel, Berlin}




\bibart{SchrSetRec}
{B. Schr\"oder}{2022}
{Set recognition of decomposable graphs and steps
towards their reconstruction}
{Abh. Math. Semin. Univ. Hambg.}{92}{1--25}

\url{https://doi.org/10.1007/s12188-021-00252-0}



\end{thebibliography}
\end{document}